\documentclass[12pt,twoside,reqno,openany]{amsart}
\pdfoutput=1
\usepackage{amsbsy,amscd,amsfonts,amsmath,amssymb,amsthm,color,
fancybox,fancyhdr,footmisc,graphics,graphicx,ifthen,mathrsfs,
multicol,pdfpages,rotating,times,wasysym}

\usepackage[dvipsnames,svgnames,x11names]{xcolor}

\usepackage[all]{xy}
\usepackage[utf8]{inputenc}
\usepackage[T1]{fontenc}
\usepackage{mathtools}
\newtagform{EngelLie}[\scriptstyle]{$}{$}
\makeatletter\newcommand{\leqnomode}{\tagsleft@true}
\newcommand{\reqnomode}{\tagsleft@false}\makeatother

\newtheorem{Theorem}[equation]{Theorem}

\newtheorem{Proposition}[equation]{Proposition}

\newtheorem{Lemma}[equation]{Lemma}

\newtheorem{Corollary}[equation]{Corollary}

\theoremstyle{definition}

\newtheorem{Definition}[equation]{Definition}

\definecolor{blue}{cmyk}{1.,1.,0.,0.63}
\definecolor{red}{cmyk}{0.,1.,1.,0.63}
\definecolor{green}{cmyk}{1.,0.,1.,0.63}
\definecolor{black}{cmyk}{1.,1.,1.,1.}

\makeatletter
\renewcommand{\@fnsymbol}[1]
{\ensuremath{\ifcase#1\or $*$\or $**$\or $***$\or $****$\or $*****$
\else\@ctrerr\fi}}
\makeatother

\numberwithin{equation}{section}

\newcommand{\style}[1]{\text{\footnotesize{\sf #1}}}

\renewcommand{\dim}{\style{dim}}

\renewcommand{\exp}{\style{exp}}

\renewcommand{\Im}{\style{Im}}

\renewcommand{\ker}{\style{ker}}

\renewcommand{\lim}{\style{lim}}

\renewcommand{\mod}{\style{mod}}

\newcommand{\vf}{\vfill\end{document}}

\usepackage{graphicx,amsmath}
\newtheorem{Summary}[equation]{Summary}
\newcounter{para}

\allowdisplaybreaks

\title[Normal Forms under Fibre-Preserving Maps]{Normal Forms of second order Ordinary Differential Equations $y_{xx}=J(x,y,y_{x})$ under Fibre-Preserving Maps}
\author{Wei Guo {\sc Foo},\quad Julien {\sc Heyd},\quad Jo\"el {\sc Merker}}

\begin{document}
\begin{abstract}
We study the equivalence problem of classifying second order ordinary differential equations $y_{xx}=J(x,y,y_{x})$ modulo fibre-preserving point transformations $x\longmapsto \varphi(x)$, $y\longmapsto \psi(x,y)$ by using Moser's method of normal forms. We first compute a basis of the Lie algebra $\frak{g}_{\{y_{xx}=0\}}$ of fibre-preserving symmetries of $y_{xx}=0$. In the formal theory of Moser's method, this Lie algebra is used to give an explicit description of the set of normal forms $\mathcal{N}$, and we show that the set is an ideal in the space of formal power series. We  then show the existence of the normal forms by studying flows of suitable vector fields with appropriate corrections by the Cauchy-Kovalevskaya theorem. As an application, we show how normal forms can be used to prove that the identical vanishing of Hsu-Kamran primary invariants directly imply that the second order differential equation is fibre-preserving point equivalent to $y_{xx}=0$.
\end{abstract}

\maketitle


\section{Introduction}

The classification of second order ordinary differential equations
$y_{xx} = J(x,y,y_x)$ under point transformations was solved by Lie
and two of his students: Tresse and Cartan, who used complementary
methods. Throughout we adopt the notation $p:=y_{x}$. 

Cartan's approach, inspired by the works of his master Darboux, was
the known {\sl method of moving frames}. Further, in 1932, based on an
observation of Segre, Cartan remarked that Lie's classification of
2\textsuperscript{nd} order ordinary differential equations can be
carried over, {\sl mutatis mutandis}, to the classification of
3-dimensional Levi non-degenerate real hypersurfaces $M^3 \subset
\mathbb{C}^{2}$ \cite{Cartan-1}.  In 1974,
Chern-Moser~{\cite{Chern-Moser-1974}} constructed a Cartan connection
associated to CR equivalences of $(2n+1)$-dimensional Levi
nondegenerate CR manifolds. More recently, this method was used by
Nurowski-Sparling \cite{Nurowski-Sparling-2003} to define an analogue
of Fefferman metric for second order ordinary differential
equations. Also interestingly, Godlinski and Nurowski
\cite{Godlinski-2008,
Godlinski-Nurowski-2009} applied the same method to solve the
equivalence problem of classifying 3\textsuperscript{rd} 
order ordinary differential
equations modulo contact,
point, and fiber-preserving transformations,
with deep connections to
Einstein-Weyl geometries.

What is probably less known is that the theory of normal forms for
second order ordinary differential equations $y_{xx} = J(x,y,p)$
in the sense of Moser was already done by Tresse in his prized  
thesis~{\cite{Tresse-1894, Tresse-1896}}. Under point transformations, he obtained:
\begin{equation*}
\begin{aligned}
y_{xx} &=\sum_{i,j,k=0}J_{i,j,k}\frac{x^{i}y^{j}p^{k}}{i!j!k!}\\
&=
J_{0,0,4}\frac{p^{4}}{4!}+J_{2,0,2}\frac{x^{2}p^{2}}{2!2!}
+
\sum_{i+j+k\geqslant 1}J_{i,j,4+k}\frac{x^{i}y^{j}p^{k}}{i!j!k!}
+
\sum_{i+j+k\geqslant 1}J_{2+i,j,2+k}\frac{x^{2+i}y^{j}p^{2+k}}{(2+i)!j!(2+k)!},
\end{aligned}
\end{equation*}
where $J_{0,0,4}$ and $J_{2,0,2}$  correspond respectively to the following relative invariants at the origin:
\begin{equation*}
\begin{aligned}
Q_{1} &:= J_{pppp},\\
Q_{2} &= \widehat{\sf D}_{x}^{2}J_{pp}-4\widehat{\sf D}_{x}J_{yp}-6J_{p}J_{pp}\widehat{\sf D}_{x}J_{pp}-3J_{y}J_{pp}+4J_{p}J_{yp}.
\end{aligned}
\end{equation*}
Here the vector field $\widehat{\sf D}_{x}$ is some modification of the total differential operator ${\sf D}_{x}$. Recently in 2017, Ottazzi-Schmalz \cite{Ottazzi-Schmalz-2017} solved  completely the problem using geometric methods of Chern-Moser \cite{Chern-Moser-1974} in the context of para-CR structures.

Inspired by Poincar\'e's works on dynamical systems and celestial
mechanics, Moser in~{\cite{Chern-Moser-1974}} handled elegantly CR
equivalences thanks to his method of {\sl normal forms}. Since then,
Moser's method has generated considerable interest amongst
mathematicians working in several complex variables and/or in CR
geometry. Usually, the process is divided into two steps:
firstly, to exhibit a {\sl formal} normalisation; and secondly, to
perform a {\sl convergent} normalisation, which is more delicate. 

Using these tools, we study the equivalence problem of classifying second order ordinary differential equations modulo {\sl fibre-preserving point transformations}, which are of the form $x\longmapsto \varphi(x)$, $y\longmapsto \psi(x,y)$. The formal normalisation consists of bringing $y_{xx}=J(x,y,p)$ to normal forms $y_{xx}=K(x,y,p)$ using  formal power series $\varphi(x)$, $\psi(x,y)$ which are not necessarily convergent. Our first result is the following 
\begin{Theorem}[stated as Theorem \ref{thm-nf}]
There exists a unique choice of formal power series $(f(x),g(x,y))$ with:
\begin{equation*}
\begin{aligned}
f(0) &= f'(0) =0,\\
g(0,0) &= g_{x}(0,0)=g_{y}(0,0)=g_{xy}(0,0)=0,
\end{aligned}
\end{equation*}
such that the map $\varphi(x)=x+f(x)$, $\psi(x,y)=y+g(x,y)$, brings $y_{xx}=J(x,y,p)$ to a normal form:
\[
y_{xx}=\sum_{i,j=0}^{\infty}K_{i,j}(y)x^{i}p^{j},
\]
where the $K_{i,j}(y)$ satisfy:
\begin{enumerate}
\item $K_{i,0}(y)\equiv 0$ for  all $i\geqslant 0$;
\item $K_{0,1}(y)\equiv 0$;
\item $K_{i,1}(0)= 0$ for all $i\geqslant 1$;
\item $K_{0,2}(y)\equiv 0$. 
\end{enumerate} 
\end{Theorem}

Here, we argue that in the real-analytic category, based on a set
of principles and guidelines that make full advantage of the formal
theory, the convergence problem can be solved in a natural way. Our construction is based on the study of flows of suitable vector fields, with appropriate corrections via the Cauchy-Kovalevskaya theorem. These give our second and the main result, which is the following 

\begin{Theorem}[stated as Theorem \ref{TE}]
There exists a convergent (analytic) fibre-preserving point
transformation that sends an analytic second order ordinary
differential equation $y_{xx}=J(x,y,p)$ to an analytic normal form 
$y_{xx}=K(x,y,p)$ satisfying (1--4).
\end{Theorem}

The procedure can be found in Summary \ref{SE}. As a result, the normal form has the following expansion:
\begin{equation*}
\begin{aligned}
y_{xx}&=K_{0,0,3}\frac{p^{3}}{3!}+
K_{1,0,2}\frac{xp^{2}}{2!}
+
K_{1,1,1}xyp
\\
&\hspace{0.5cm}
+\sum_{i+j+k\geqslant 1}
K_{i,j,3+k}\frac{x^{i}y^{j}p^{3+k}}{i!j!(3+k)!}
+
\sum_{i+j+k\geqslant 1}
K_{1+i,j,2+k}\frac{x^{1+i}y^{j}p^{2+k}}{(1+i)!j!(2+k)!}\\
&\hspace{0.5cm}
+
\sum_{i+j+k\geqslant 1}
K_{1+i,1+j,1+k}\frac{x^{1+i}y^{1+j}p^{1+k}}{(1+i)!(1+j)!(1+k)!},
\end{aligned}
\end{equation*}
where the coefficients $K_{0,0,3}$, $K_{1,0,2}$, $K_{1,1,1}$ correspond respectively to the values of the relative differential invariants $I_{1}(K)$, $I_{2}(K)$, $I_{3}(K)$ (see equation \eqref{HK-rI}) at the origin, which are the only primary ones discovered by Hsu-Kamran \cite{Hsu-Kamran-1989}.

The theory of normal forms is  one of the important tools in the
theory of differential invariants. The Taylor coefficients of a normal
form capture the values of all invariants at the origin, hence they
effectively characterise the geometry of the problem. In several
articles, Olver and his collaborators gave what they call the {\sl
recurrence formulae} to express the higher order differential
invariants as derivatives of certain primary invariants \cite{Olver-Pohjanpelto-2008, Olver-2018, Olver-Valiquette-2018}. This is not
the focus of the paper for now, but we will show in our future
publications how normal forms 
can be used to find all the homogeneous models
classified by Hsu-Kamran in~{\cite{Hsu-Kamran-1989}}.
\medskip 

\noindent \textbf{Acknowledgements:} The authors would like to thank
Professors Pawe\l{} Nurowski (Centrum Fizyki Teoretycznej, Polska
Akademia Nauk) and Chin-Yu Hsiao (Institute of Mathematics, Academia
Sinica) for numerous fruitful exchanges. This research was supported
in part by the Polish National Science Centre (NCN) via the grant
number 2018/29/B/ST1/02583, and by the Norwegian Financial Mechanism
2014–2021 via the project registration number 2019/34/H/ST1/00636. Wei Guo Foo is partially supported by Taiwan Ministry of Science and Technology projects  108-2115-M-001-012-MY5, 109-2923-M-001-010-MY4.


\section{Summary of the result of Hsu-Kamran}

The problem of our study is the classification of second order ordinary differential equations
\[
Y_{XX}=J(X,Y,Y_{X}),
\]
up to fibre-preserving maps $\Phi$:
\begin{equation}
\begin{aligned}
X &= \varphi(x),\\
Y &= \psi(x,y).
\end{aligned}
\end{equation}
This problem was solved by {\sc Hsu-Kamran} (for example \cite{Hsu-Kamran-1989}) using method of moving frames, and we briefly summarise the results leading up to the existence of three relative invariants. 

Let $J^{2}(\mathbb{R},\mathbb{R})$ be the second order Jet space with independent coordinates $(x,y,y_{x},y_{xx})$. For simplicity, we adopt the notation $p:=y_{x}$. The total differential operator ${\sf D}_{x}$ along the $x$-axis is 
\[
{\sf D}_{x} := \frac{\partial}{\partial x}+
p\frac{\partial}{\partial y}+
y_{xx}\frac{\partial}{\partial p}
+\sum_{j=3}^{+\infty}y_{x^{(k)}}\frac{\partial}{\partial y_{x^{(k-1)}}}.
\]
The fibre preserving maps $\Phi$ can be prolonged to second order 
\begin{equation*}
\begin{aligned}
\Phi^{2}:\ J^{2}(\mathbb{R},\mathbb{R}) &\longmapsto J^{2}(\mathbb{R},\mathbb{R})\\
(x,y,p,y_{xx}) &\longmapsto (X,Y,P,Y_{XX})
\end{aligned}
\end{equation*}
using the recursive formula well-known to {\sc Sophus Lie}:
\begin{equation}\label{ptr-02}
\begin{aligned}
X &= \varphi(x),\\
Y &= \psi(x,y),\\
P &=\frac{{\sf D}_{x}Y}{{\sf D}_{x}X}= \frac{\psi_{x}+\psi_{y}p}{\varphi_{x}},\\
Y_{XX} &= \frac{{\sf D}_{x}P}{{\sf D}_{x}X}=
\frac{\varphi_{x}\psi_{yy}p^{2}+2\varphi_{x}\psi_{xy}p+\varphi_{x}\psi_{xx}+(-\psi_{x}-\psi_{y}p)\varphi_{xx}+\psi_{y}\varphi_{x}y_{xx}}{\varphi_{x}^{3}}.
\end{aligned}
\end{equation}

A given second order differential equation $y_{xx}=K(x,y,p)$ defines a 3-dimensional submanifold $M\subset J^{2}(\mathbb{R},\mathbb{R})$. The cotangent bundle $T^{*}M$ is generated by the following $1$-forms:
\begin{equation}
\begin{aligned}
\omega^{1} &:= dx,\\
\omega^{2} &:= dy-p\ dx,\\
\omega^{3} &:= dp-K(x,y,p)\ dx.
\end{aligned}
\end{equation}
Now, suppose that $\Phi^{2}$ sends $M$ to $M'$ given by $Y_{XX}=J(X,Y,P)$. A simple substitution by \eqref{ptr-02} shows that the pull-back $\Phi^{2}$ transfers $\{\omega^{1},\omega^{2},\omega^{3}\}$ to $\{\Omega^{1}:= dX,\ \Omega^{2} := dY-P\ dX,\ \Omega^{3}:= dP-J(X,Y,P)\ dX\}$ via the initial ${\sf G}$-structure:
\begin{equation}
\begin{aligned}
\left(
\begin{matrix}
\Omega^{1}\\
\Omega^{2}\\
\Omega^{3}
\end{matrix}
\right)
=
\left(
\begin{matrix}
{\sf A} & 0 & 0\\
0 & {\sf B} & 0\\
0 & {\sf BC} & {\sf B}/{\sf A}
\end{matrix}
\right)
\left(
\begin{matrix}
\omega^{1}\\
\omega^{2}\\
\omega^{3}
\end{matrix}
\right).
\end{aligned}
\end{equation}

By applying {\sc Cartan}'s method of moving frames, {\sc Hsu-Kamran} obtained a 6-dimensional principal bundle ${\sf P}\rightarrow M$ whose cotangent bundle is generated by six invariant $1$-forms $\{\omega^{1},\omega^{2},\omega^{3},\alpha,\beta,\gamma\}$ satisfying the structure equations:
\begin{equation}
\begin{aligned}
d\omega^{1} &= \alpha\wedge\omega^{1},\\
d\omega^{2} &= \beta\wedge\omega^{2}+\omega^{1}\wedge\omega^{3},\\
d\omega^{3} &= \gamma\wedge\omega^{2}+(\beta-\alpha)\wedge\omega^{3},\\
d\alpha &= -2\gamma\wedge\omega^{1},\\
d\beta &= \omega^{1}\wedge\gamma+ I_{2}\ \omega^{1}\wedge\omega^{2}-I_{1}\ \omega^{3}\wedge\omega^{2},\\
d\gamma &= \gamma\wedge\alpha + I_{3}\ \omega^{1}\wedge\omega^{2}+I_{2}\ \omega^{1}\wedge\omega^{3}.
\end{aligned}
\end{equation}
Parametrically, there are 3 basic relative invariants: 
\begin{equation}\label{HK-rI}
\begin{aligned}
I_{1}(K) &= -\frac{{\sf A}}{2{\sf B}^{2}}K_{ppp},\\
I_{2}(K) &= \frac{1}{2{\sf AB}}({\sf D}_{x}K_{pp}-K_{py}),\\
I_{3}(K) &= -{\sf C}I_{2}+\frac{1}{2{\sf A}^{2}{\sf B}}({\sf D}_{x}K_{py}+K_{pp}K_{y}-K_{py}K_{p}-2K_{yy}),
\end{aligned}
\end{equation}
while the invariant $1$-forms $\alpha$, $\beta$, $\gamma$ are:
\begin{equation}
\begin{aligned}
\alpha &= \frac{d{\sf A}}{\sf A}-\bigg(2{\sf C}+\frac{K_{p}}{\sf A}\bigg)\omega^{1},\\
\beta &= \frac{d{\sf B}}{\sf B}-{\sf C}\omega^{1}+\frac{K_{pp}}{2{\sf B}}\omega^{2},\\
\gamma &= d{\sf C}+{\sf C}\frac{d{\sf A}}{\sf A}+
\bigg(\frac{K_{y}}{{\sf A}^{2}}-
\frac{{\sf C}K_{p}}{\sf A}
-{\sf C}^{2}\bigg)\omega^{1}
+
\bigg(
\frac{K_{py}}{2{\sf AB}}-\frac{{\sf C}K_{pp}}{2{\sf B}}\bigg)\omega^{2}
+
\frac{K_{pp}}{2{\sf B}}\omega^{3}.
\end{aligned}
\end{equation}

\section{Normal Forms: Preliminaries}

\subsection{} Let ${\sf G}$ be the infinite-dimensional Lie pseudo-group of all fibre-preserving point transformations, and let $\frak{g}$ be its Lie algebra. An element ${\sf X}\in\frak{g}$ is a vector field on $\mathbb{R}^{2}$ of the form:
\[
{\sf X} = f(x)\frac{\partial}{\partial x}+g(x,y)\frac{\partial}{\partial y}.
\]
The 1-parameter flow $\exp(t{\sf X})(x,y)$ of ${\sf X}$ is obtained by solving the two ordinary differential equations:
\begin{equation}
\begin{aligned}
\gamma_{1}'(t) &= f(\gamma_{1}(t)),\\
\gamma_{2}'(t) &= g(\gamma_{1}(t),\gamma_{2}(t)),
\end{aligned}
\end{equation}
with the initial conditions $\gamma_{1}(0)=x$, $\gamma_{2}(0)=y$. Clearly, the resulting $1$-parameter family of fibre-preserving point transformations has the  expansion in $t$:
\begin{equation}
\begin{aligned}
X &= x + f(x)t+ O(t^{2}),\\
Y &= y+g(x,y)t + O(t^{2}). 
\end{aligned}
\end{equation}

The vector field can be prolonged to second order
\[
{\sf X}^{(2)}= 
{\sf X}+{\sf X}_{1}(x,y,p)\frac{\partial}{\partial p}+
{\sf X}_{2}(x,y,p,y_{xx})\frac{\partial}{\partial y_{xx}},
\]
via the well-known prolongation formula:
\begin{equation}
\begin{aligned}
{\sf X}_{1} &= {\sf D}_{x}\big(g(x,y)-f(x)y\big)-f(x)y_{xx},\\
{\sf X}_{2} &= {\sf D}_{x}^{2}\big(g(x,y)-f(x)y\big)-f(x)y_{xxx}.
\end{aligned}
\end{equation}

\subsection{} The vector field ${\sf X}$ is a fibre-preserving infinitesimal symmetry of a second order ordinary differential equation $y_{xx}=K(x,y,p)$ if and only if its second prolongation ${\sf X}^{(2)}$ is tangent to the manifold $M\subset J^{2}(\mathbb{R},\mathbb{R})$, namely:
\[
{\sf X}^{(2)}(y_{xx}-K(x,y,p))\big|_{y_{xx}=K(x,y,p)}\equiv 0.
\]
The set of such vector fields is denoted by $\frak{g}_{M}$, and denote its second prolongation by $\frak{g}_{M}^{(2)}$. The most important model case is when $M=\{y_{xx}=0\}$. 
\begin{Theorem}\label{yxx-sym}
Let $M=\{y_{xx}=0\}$ be a second order ordinary differential equation. Then $\dim\frak{g}_{M}=6$, and its second prolongation $\frak{g}_{M}^{(2)}$ is generated by
\begin{equation}
\begin{aligned}
{\sf v}_{1}^{(2)} &= \partial_{x},\\
{\sf v}_{2}^{(2)} &= \partial_{y},\\
{\sf v}_{3}^{(2)} &= x\partial_{y}+\partial_{p},\\
{\sf v}_{4}^{(2)} &= y\partial_{y}+p\partial_{p}+y_{xx}\partial_{y_{xx}},\\
{\sf v}_{5}^{(2)} &= x\partial_{x}-p\partial_{p}-2y_{xx}\partial_{y_{xx}},\\
{\sf v}_{6}^{(2)} &= x^{2}\partial_{x}+xy\partial_{y}-(2xp-y)\partial_{p}-3xy_{xx}\partial_{y_{xx}}.
\end{aligned}
\end{equation}
\end{Theorem}
\begin{proof}
By expanding the left-hand side of 
\[
{\sf X}^{(2)}(y_{xx})\big|_{y_{xx}=0}\equiv 0,
\]
we obtain
\begin{equation}
\begin{aligned}
g_{xx}+(2g_{xy}-f_{xx})p+g_{yy}p^{2}=0.
\end{aligned}
\end{equation}
Solving for $f(x)$, $g(x,y)$ the following system of partial differential equations:
\begin{equation*}
\begin{aligned}
g_{xx} &= 0,\\
2g_{xy} &= f_{xx},\\
g_{yy} &= 0,
\end{aligned}
\end{equation*}
the vector field ${\sf X}$ can be written as a linear combination of the following six vector fields:
\begin{equation}
\begin{aligned}
{\sf v}_{1} &= \partial_{x},\\
{\sf v}_{2} &= \partial_{y},\\
{\sf v}_{3} &= x\partial_{y},\\
{\sf v}_{4} &= y\partial_{y},\\
{\sf v}_{5} &= x\partial_{x},\\
{\sf v}_{6} &= x^{2}\partial_{x}+xy\partial_{y}.
\end{aligned}
\end{equation}
Prolonging all of them to second order finishes the proof. 
\end{proof}

\subsection{} The equation
\[
g_{xx}+(2g_{xy}-f_{xx})p+g_{yy}p^{2}=0,
\]
is sometimes called the defining function of $\frak{g}_{\{y_{xx}=0\}}$ in $\frak{g}$, and it will also be used to define our normal form. In fact, let $\mathcal{F}_{x}$ denote the space of all formal power series in $x$, let $\mathcal{F}_{x,y}$ denote the space of all formal power series in $(x,y)$, and let $\mathcal{F}_{x,y,p}$ denote the space of all formal power series in $(x,y,p)$. Introduce the homological operator:
\begin{equation}
\begin{aligned}
L:\hspace{0.75cm} \mathcal{F}_{x}\times\mathcal{F}_{x,y} &\longrightarrow \mathcal{F}_{x,y,p}\\
(f(x),g(x,y)) &\longmapsto L(f,g):=g_{xx}+(2g_{xy}-f_{xx})p+g_{yy}p^{2}.
\end{aligned}
\end{equation}
We seek a subspace $\mathcal{N}\subset \mathcal{F}_{x,y,p}$ which we call the space of normal forms, satisfying:
\[
\mathcal{F}_{x,y,p}=\Im(L)\oplus \mathcal{N},
\]
where $\mathcal{N}$ consists of representatives $K$ of classes $[K]$ in $\mathcal{F}_{x,y,p}/\Im(L)$ whose image part has been completely normalised by $L(f,g)$ for some $(f,g)\in\mathcal{F}_{x}\times\mathcal{F}_{x,y}$ (or in other words, has been completely absorbed into $\Im(L)$).

 A question is: \textit{how can we find an explicit description of $\mathcal{N}$?} A clue is to look at the kernel of $L$. To say that there are two different ways to bring $y_{xx}=J(x,y,p)$ to a normal form is to say that there exist two different $(f,g)$ and $(\widehat{f}, \widehat{g})$ in $\mathcal{F}_{x}\times\mathcal{F}_{x,y}$ such that 
\begin{equation}
\begin{aligned}
J(x,y,p) &= K(x,y,p)+L(f,g),\\
J(x,y,p) &= K(x,y,p)+L(\widehat{f},\widehat{g}),
\end{aligned}
\end{equation}
with $K$ in $\mathcal{N}$. Then 
\[
L(\widehat{f}-f,\widehat{g}-g)=0,
\]
or 
\[
(\widehat{f}-f,\widehat{g}-g)\in \ker\ L.
\]
Therefore, the choices of normalisations are unique up to {\sl elements in $\ker\ L$}. Fortunately, we have information about this indeterminacy because from Theorem \ref{yxx-sym},  $(f(x),g(x,y))\in \ker\ L$ if and only if the corresponding vector field ${\sf X}=f(x)\partial_{x}+g(x,y)\partial_{y}$ is a fibre-preserving infinitesimal symmetry of $y_{xx}=0$, and we have an explicit basis of this Lie algebra $\{{\sf v}_{1},\ {\sf v}_{2},\ {\sf v}_{3},\ {\sf v}_{4},\ {\sf v}_{5},\ {\sf v}_{6}\}$. 

This piece of information can help us find the subspace $\mathcal{F}'\subset \mathcal{F}_{x}\times\mathcal{F}_{x,y}$ of the source space on which $L$ is \textit{injective}.  Then any element $(f,g)\in\mathcal{F}'$ coming from this subspace will give the {\sl unique} normalisation, hence giving us a precise description of a normal form. To put this idea into action, we first have to find out what $\mathcal{F}'$ is. We expand ${\sf X}$ in terms of power series:
\[
{\sf X} = 
\big(\alpha_{0}+\alpha_{1}x+\alpha_{2}x^{2}+\alpha_{3}x^{3}+O_{x}(4)\big)\frac{\partial}{\partial x}
+
\big(\beta_{0,0}+\beta_{1,0}x + \beta_{0,1}y + \beta_{2,0}x^{2}+\beta_{1,1}xy+\beta_{0,2}y^{2}+O_{x,y}(3)\big)\frac{\partial}{\partial y}.
\]
Modulo the vector space generated by the first $5$ vector fields ${\sf v}_{1}$, ${\sf v}_{2}$, ${\sf v}_{3}$, ${\sf v}_{4}$, ${\sf v}_{5}$:
\[
{\sf X} \equiv 
\big(\alpha_{2}x^{2}+\alpha_{3}x^{3}+O_{x}(4)\big)\frac{\partial}{\partial x}
+
\big( \beta_{2,0}x^{2}+\beta_{1,1}xy+\beta_{0,2}y^{2}+O_{x,y}(3)\big)\frac{\partial}{\partial y}\ \mod\ \langle {\sf v}_{1},\ {\sf v}_{2},\ {\sf v}_{3}, {\sf v}_{4}, {\sf v}_{5}\rangle,
\]
which can be written either as:
\begin{equation}
\begin{aligned}
{\sf X} &\equiv 
\big((\alpha_{2}-\beta_{1,1})x^{2}+\alpha_{3}x^{3}+O_{x}(4)\big)\frac{\partial}{\partial x}
+
\big( \beta_{2,0}x^{2}+\beta_{0,2}y^{2}+O_{x,y}(3)\big)\frac{\partial}{\partial y}\\
&\hspace{0.5cm}
+\beta_{1,1}{\sf v}_{6}\quad
 \mod\ \langle {\sf v}_{1},{\sf v}_{2},{\sf v}_{3}, {\sf v}_{4}, {\sf v}_{5}\rangle,
\end{aligned}
\end{equation}
or
\begin{equation}
\begin{aligned}
{\sf X} &\equiv 
\big(\alpha_{3}x^{3}+O_{x}(4)\big)\frac{\partial}{\partial x}
+
\big( \beta_{2,0}x^{2}+(\beta_{1,1}-\alpha_{2})xy+\beta_{0,2}y^{2}+O_{x,y}(3)\big)\frac{\partial}{\partial y}\\
&\hspace{0.5cm}
+\alpha_{2}{\sf v}_{6}\quad
 \mod\ \langle {\sf v}_{1}, {\sf v}_{2},{\sf v}_{3}, {\sf v}_{4}, {\sf v}_{5}\rangle.
\end{aligned}
\end{equation}
The two choices are valid, and we will choose the first one for our construction. Hence 
\begin{equation*}
\begin{aligned}
{\sf X} &\equiv 
\big((\alpha_{2}-\beta_{1,1})x^{2}+\alpha_{3}x^{3}+O_{x}(4)\big)\frac{\partial}{\partial x}
+
\big( \beta_{2,0}x^{2}+\beta_{0,2}y^{2}+O_{x,y}(3)\big)\frac{\partial}{\partial y}\quad \mod\ \frak{g}_{\{y_{xx}=0\}}.
\end{aligned}
\end{equation*}
Thus we have the following
\begin{Proposition}\label{F-prime}
The subspace $\mathcal{F}'\subset \mathcal{F}_{x}\times \mathcal{F}_{x,y}$ on which the homological operator $L$ is injective is given by the set of all tuples $(f(x),g(x,y))$ satisfying the following conditions:
\begin{enumerate}
\item  $f(x)$ starts with $x^{2}$,
\item  $g(x,y)$  starts with homogeneous terms of order $2$, with $g_{xy}(0,0)=0$. \qed
\end{enumerate}
\end{Proposition}
Using $\mathcal{F}'$, we may now proceed with the computation of $\mathcal{N}$.  

\section{Explicit formal computations of normal forms $\mathcal{N}$}

We expand $(f(x),g(x,y))\in\mathcal{F}'$ in terms of their power series:
\begin{equation}
\begin{aligned}
f(x) &=\sum_{j=2}^{\infty}f_{j}x^{2},\\
g(x,y) &= \sum_{j=0}^{\infty}g_{j}(y)x^{j},
\end{aligned}
\end{equation}
with, due to Proposition \ref{F-prime}:
\[
g_{0}(0)=g_{0}'(0)=0,\qquad 
g_{1}(0)=g_{1}'(0)=0.
\]
We write $K(x,y,p)$ and $J(x,y,p)$ as 
\begin{equation}
\begin{aligned}
K(x,y,p)=\sum_{i,j=0}^{\infty} K_{i,j}(y)\ x^{i}p^{j},\\
J(x,y,p)=\sum_{i,j=0}^{\infty} J_{i,j}(y)\ x^{i}p^{j}.
\end{aligned}
\end{equation}
The homological equation
\[
K(x,y,p)=J(x,y,p)+L(f,g)
\]
becomes
\begin{equation*}
\begin{aligned}
\sum_{i,j=0}^{\infty} K_{i,j}(y)\ x^{i}p^{j}
&=
\sum_{i,j=0}^{\infty} J_{i,j}(y)\ x^{i}p^{j}
+
\sum_{i=0}^{\infty}(i+1)(i+2)g_{i+2}(y)x^{i}\\
&\hspace{0.5cm}
+
2\bigg(
\sum_{i=0}^{\infty}(i+1)g_{i+1}'(y)x^{i}
-
\sum_{i=0}^{\infty}(i+1)(i+2)f_{i+2}x^{i}\bigg)p
+
\bigg(\sum_{i=0}^{\infty}g_{i}''(y)x^{i}\bigg)p^{2}.
\end{aligned}
\end{equation*}
It is therefore clear that only terms containing $p$ of degrees $0$, $1$, $2$ will be subjected to change during the formal computations. For the terms not containing $p$, with any $i\geqslant 0$:
\begin{equation}\label{A}
K_{i,0}(y) = J_{i,0}(y)+(i+1)(i+2)g_{i+2}(y).
\end{equation}
For the terms containing $p$, with any $i\geqslant 0$,
\begin{equation}\label{B}
K_{i,1}(y)=J_{i,1}(y)+2(i+1)g_{i+1}'(y)-2(i+1)(i+2)f_{i+2}.
\end{equation}
Finally for terms containing $p^{2}$, with any $i\geqslant 0$,
\begin{equation}\label{C}
K_{i,2}(y)=J_{i,2}(y)+g_{i}''(y).
\end{equation}
We are ready to find an explicit description of $\mathcal{N}$. 
\begin{Theorem}\label{thm-nf}
There exists a unique choice of $(f(x),g(x,y))\in\mathcal{F}'$ such that the $K_{i,j}(y)$ in 
\[
y_{xx}=\sum_{i,j=0}^{\infty}K_{i,j}(y)x^{i}p^{j},
\]
satisfy:
\begin{enumerate}
\item $K_{i,0}(y)\equiv 0$ for  all $i\geqslant 0$;
\item $K_{0,1}(y)\equiv 0$;
\item $K_{i,1}(0)= 0$ for all $i\geqslant 1$;
\item $K_{0,2}(y)\equiv 0$. 
\end{enumerate} 
\end{Theorem}
\begin{proof}
The first condition can be obtained from equation \eqref{A} by choosing for $i\geqslant 0$:
\begin{equation}\label{sol-1}
g_{i+2}(y)=-\frac{1}{(i+1)(i+2)}J_{i,0}(y).
\end{equation}
To get (2), set $i=0$ in equation \eqref{B}:
\[
4f_{2}+2g_{1}'(y)=-J_{0,1}(y).
\]
When $y=0$, the hypothesis $g_{1}'(0)=0$ gives the value for $f_{2}$:
\begin{equation}\label{sol-2-f}
f_{2}=-\tfrac{1}{4}J_{0,1}(0),
\end{equation}
from which the value for $g_{1}'(y)$ can be deduced:
\begin{equation}\label{sol-2-g}
g_{1}'(y)=\tfrac{1}{2}(-J_{0,1}(y)+J_{0,1}(0)).
\end{equation}
To get (3), again from equation \eqref{B} for $i\geqslant 1$,
\begin{equation*}
\begin{aligned}
K_{i,1}(0) &= 
J_{i,1}(0)+2(i+1)g_{i+1}'(0)-2(i+1)(i+2)f_{i+2}\\
&=
J_{i,1}(0)+2(i+1)\bigg(\frac{-1}{i(i+1)}J_{i-1,0}'(0)\bigg)-2(i+1)(i+2)f_{i+2}\\
&=
J_{i,1}(0)-\tfrac{2}{i}J_{i-1}'(0)-2(i+1)(i+2)f_{i+2},
\end{aligned}
\end{equation*}
which can be normalised to zero if we choose
\begin{equation}\label{sol-3}
f_{i+2}=\tfrac{1}{(i+1)(i+2)}\big(J_{i,1}(0)-\tfrac{2}{i}J_{i-1,0}'(0)\big).
\end{equation}
Finally for (4), we can only eliminate $K_{0,2}(y)$ since the only remaining free variable is $g_{0}''(y)$. From \eqref{C},
\begin{equation}\label{sol-4}
K_{0,2}(y)=J_{0,2}(y)+g_{0}''(y),
\end{equation}
and this finishes the proof. 
\end{proof}

Consequently,  the function $K(x,y,p)\in\mathcal{N}$ is in normal form if and only if the four conditions in the theorem above are satisfied. This motivates the following Definition (or rather a Theorem--Definition):
\begin{Definition}
A second order ordinary differential equation $y_{xx}=K(x,y,p)$ is in normal form under fibre-preserving maps (i.e. $K\in\mathcal{N}$)  if and only if $K(x,y,p)$ satisfies the following set of partial differential equations:
\begin{itemize}
\item[D1.] $K(x,y,0) \equiv 0$,
\item[D2.] $K_{p}(0,y,0) \equiv 0$,
\item[D3.] $K_{p}(x,0,0)\equiv 0$,
\item[D4.] $K_{pp}(0,y,0)\equiv 0$. 
\end{itemize}
\end{Definition}
We leave it as an exercise to the readers to verify that the conditions (1), (2), (3), (4) in Theorem \ref{thm-nf} correspond exactly to these conditions D1, D2, D3, D4. This definition is useful in many ways because this formulation is in closed form and is precise. Firstly, it will be the main ingredient in the proof of the convergence of the normal form. Secondly, in Chern-Moser, it was stated {\sl without proof} that the set of normal forms generate an ideal in the ring of formal power series. Let us prove this result in the context of fibre-preserving equivalences:
\begin{Proposition}
The set of normal forms $\mathcal{N}$ is an ideal in $\mathcal{F}_{x,y,p}$. 
\end{Proposition}
\begin{proof}
It is clear that $\mathcal{N}$ forms a vector subspace of $\mathcal{F}_{x,y,p}$. It suffices to show that for any $A(x,y,p)\in \mathcal{F}_{x,y,p}$, and for any $K(x,y,p)\in\mathcal{N}$, we have $B(x,y,p):=A(x,y,p)K(x,y,p)\in\mathcal{N}$. We have to verify that the product satisfies D1--D4 conditions. 

For the condition D1:
\[
B(x,y,0)=A(x,y,0)K(x,y,0)\equiv 0. 
\]
For the conditions D2 and D3, we differentiate $B(x,y,p)$ with respect to $p$:
\[
B_{p}(x,y,p)=A_{p}(x,y,p)K(x,y,p)+A(x,y,p)K_{p}(x,y,p).
\]
Hence D2 is satisfied:
\[
B_{p}(0,y,0)=A_{p}(0,y,0)K(0,y,0)+A(0,y,0)K_{p}(0,y,0)\equiv 0,
\]
and so is D3:
\[
B_{p}(x,0,0)=A_{p}(x,0,0)K(x,0,0)+A(x,0,0)K_{p}(x,0,0)\equiv 0. 
\]
Finally for the condition D4, we differentiate $B_{p}$ with respect to $p$ once more:
\[
B_{pp}(x,y,p)=A_{pp}(x,y,p)K(x,y,p)+2A_{p}(x,y,p)K_{p}(x,y,p)+A(x,y,p)K_{pp}(x,y,p).
\]
Thus
\[
B_{pp}(0,y,0)=A_{pp}(0,y,0)K(0,y,0)+2A_{p}(0,y,0)K_{p}(0,y,0)+A(0,y,0)K_{pp}(0,y,0)\equiv 0,
\]
hence the proof is completed. 
\end{proof}

\begin{Proposition}
The second order ordinary differential equation $y_{xx}=K(x,y,p)$ is in a normal form under fibre-preserving maps if and only if it can be written as 
\begin{equation}
\begin{aligned}
y_{xx}=
\sum_{i,j,k=0}^{\infty}K_{i,j,3+k}\frac{x^{i}y^{j}p^{3+k}}{i!j!(3+k)!}
+
\sum_{i,j=0}^{\infty}K_{1+i,j,2}\frac{x^{1+i}y^{j}p^{2}}{(1+i)!j!2!}
+
\sum_{i,j=0}^{\infty}K_{1+i,1+j,1}\frac{x^{1+i}y^{1+j}p}{(1+i)!(1+j)!}.
\end{aligned}
\end{equation}
\end{Proposition}
\begin{proof}
The condition D1 implies that 
\[
K(x,y,p)=A(x,y,p)p.
\]
Differentiating both sides with respect to $p$:
\[
K_{p}(x,y,p)=A_{p}(x,y,p)p+A(x,y,p).
\]
The conditions D2 and D3 imply that 
\[
A(0,y,0)\equiv 0,\qquad A(x,0,0)\equiv 0. 
\]
The first identity gives
\[
A(x,y,p)=A_{1}(x,y,p)x+A_{2}(x,y,p)p,
\]
while the second identity implies that 
\[
A_{1}(x,0,0)x\equiv 0,
\]
whence 
\[
A_{1}(x,y,p)=A_{3}(x,y,p)y+A_{4}(x,y,p)p.
\]
Thus
\[
A(x,y,p)=A_{3}(x,y,p)xy+A_{4}(x,y,p)xp+A_{2}(x,y,p)p.
\]
The final condition D4 implies that 
\[
A_{2}(0,y,0)\equiv 0,
\]
so that 
\[
A_{2}(0,y,0)=A_{5}(x,y,p)x+A_{6}(x,y,p)p.
\]
Putting these information together, we obtain
\[
K(x,y,p)=A_{3}(x,y,p)xyp+(A_{4}(x,y,p)+A_{5}(x,y,p))xp^{2}+A_{6}(x,y,p)p^{3},
\]
which finishes the proof. 
\end{proof}

\section{Discussion on the Uniqueness of Normal Forms}

This formalism via the homological operator, very importantly, tells us what are the free parameters at each stage of the normalisation process, but it {\sl hides} the complexity of the actual computations by ignoring coefficients of lower order terms in $f(x)$ and $g(x,y)$ that have been used up in the previous stages. Since they are no longer free, they are conveniently treated as zero in this formal theory. 

To see the actual computations in full force, we look at the fibre-preserving maps that are close to the identity:
\begin{equation}
\begin{aligned}
X &= x+f(x),\\
Y &= y+g(x,y),
\end{aligned}
\end{equation}
with $f_{x}(0)=0$ and $g_{y}(0,0)=0$.  By the usual prolongation formula, the prolonged map to second order has the form:
\begin{equation}
\begin{aligned}
P &= \tfrac{1}{1+f_{x}}\big(g_{x}+(1+g_{y})p\big),\\
Y_{XX} &=\tfrac{1}{1+f_{x}}\big(g_{xx}f_{x}-g_{x}f_{xx}+\boxed{g_{xx}}\\
&\hspace{1.5cm}
+(2g_{xy}f_{x}-g_{y}f_{xx}+\boxed{2g_{xy}-f_{xx}})p
+
g_{yy}(\boxed{1}+f_{x})p^{2}\\
&\hspace{1.5cm}
+
(g_{y}f_{x}+g_{y}+f_{x}+1)y_{xx}\big).
\end{aligned}
\end{equation}
We will substitute these into $Y_{XX}=J(X,Y,P)$. First, we look at the term $J(X,Y,P)$. By a standard formula:
\begin{equation}
\begin{aligned}
J(X,Y,P) &= J\big(x+f, y+g, p+(P-p)\big)\\
&=
J(x,y,p)+
\sum_{i+j+k\geqslant 1}^{\infty}
J_{x^{i}y^{j}p^{k}}(x,y,p)
\frac{f^{i}g^{j}(P-p)^{k}}{i!j!k!}.
\end{aligned}
\end{equation}
Hence 
\begin{equation*}
\begin{aligned}
(1+f_{x})J(X,Y,P) &= 
(1+f_{x})\bigg[J(x,y,p)+
\sum_{i+j+k\geqslant 1}^{\infty}
J_{x^{i}y^{j}p^{k}}(x,y,p)
\frac{f^{i}g^{j}(P-p)^{k}}{i!j!k!}\bigg]\\
&=
J(x,y,p)+
\underbrace{
J(x,y,p)f_{x}
+
(1+f_{x})
\sum_{i+j+k\geqslant 1}^{\infty}
J_{x^{i}y^{j}p^{k}}(x,y,p)
\frac{f^{i}g^{j}(P-p)^{k}}{i!j!k!}}_{\text{Non-linear terms.}}.
\end{aligned}
\end{equation*}
Next, we combine this with the expression of $Y_{XX}$, we see that 
\begin{equation*}
\begin{aligned}
y_{xx}&=
\tfrac{1}{1+f_{x}+g_{y}+g_{y}f_{x}}
\bigg[
J(x,y,p)
-
\boxed{\big(g_{xx}+(2g_{xy}-f_{xx})p+g_{yy}p^{2}\big)}\\
&\hspace{3cm}
-g_{xx}f_{x}-g_{x}f_{xx}
-
(2g_{xy}f_{x}-g_{y}f_{xx})p-g_{yy}f_{x}p^{2}\\
&\hspace{3cm}
+J(x,y,p)f_{x}
+
(1+f_{x})
\sum_{i+j+k\geqslant 1}^{\infty}
J_{x^{i}y^{j}p^{k}}(x,y,p)
\frac{f^{i}g^{j}(P-p)^{k}}{i!j!k!}\bigg].
\end{aligned}
\end{equation*}
We expand the factor in terms of power series:
\[
\frac{1}{1+f_{x}+g_{y}+g_{y}f_{x}}
=
1+\sum_{k=1}^{\infty}\frac{(-f_{x}-g_{y}-g_{y}f_{x})^{k}}{k!},
\]
and group $J(x,y,p)$ together with the boxed terms, we obtain the {\sl fundamental identity}:
\begin{equation}\label{FE}
\begin{aligned}
y_{xx} &=
J(x,y,p)
-
\boxed{\big(g_{xx}+(2g_{xy}-f_{xx})p+g_{yy}p^{2}\big)}\\
&\hspace{0.5cm}
+\bigg(J(x,y,p)
-
\big(g_{xx}+(2g_{xy}-f_{xx})p+g_{yy}p^{2}\big)\bigg)
\sum_{k=1}^{\infty}\frac{(-f_{x}-g_{y}-g_{y}f_{x})^{k}}{k!}\\
&\hspace{0.5cm}
+
\bigg(
-g_{xx}f_{x}-g_{x}f_{xx}
-
(2g_{xy}f_{x}-g_{y}f_{xx})p-g_{yy}f_{x}p^{2}+J(x,y,p)f_{x}\\
&\hspace{1cm}
+
(1+f_{x})
\sum_{i+j+k\geqslant 1}^{\infty}
J_{x^{i}y^{j}p^{k}}(x,y,p)
\frac{f^{i}g^{j}(P-p)^{k}}{i!j!k!}
\bigg)\sum_{k=0}^{\infty}\frac{(-f_{x}-g_{y}-g_{y}f_{x})^{k}}{k!}\\
&=:
J(x,y,p)-L(f,g)+R(f,g,J)\\
&=:
K(x,y,p).
\end{aligned}
\end{equation}

We will stratify the right-hand side of the equation \eqref{FE} above by assigning suitable weights to the variables $x$, $y$, $p$. From the list of symmetries of $y_{xx}=0$, we observe that 
\[
{\sf v}_{4}^{(2)}+2{\sf v}_{5}^{(2)}=x\partial_{x}+2y\partial_{y}+p\partial_{p}+0\cdot \partial_{y_{xx}}.
\]
Computing the flow of this vector field, the second order ordinary differential equation $y_{xx}=0$ is invariant under the scaling map:
\[
x \longmapsto tx,\quad y\longmapsto t^{2}y,\quad p\longmapsto tp\eqno
({\scriptstyle t>0}).
\]
Thus it would be convenient to assign the following weights: 
\[
[x]=1,\qquad 
[y]=2,\qquad 
[p]=1.
\]

Based on the weight assignments,  a function $F_{k}(x,y,p)$ is said to be {\sl semi-homogeneous of order $k$} if for all $t>0$, 
\[
F_{k}(tx,t^{2}y, tp)=t^{k}F_{k}(x,y,p).
\]
The functions $K(x,y,p)$, $J(x,y,p)$ may then be written as a sum of semi-homogeneous terms:
\begin{equation*}
\begin{aligned}
K(x,y,p)&=\sum_{\alpha=0}^{+\infty}K_{\alpha}(x,y,p),\\
J(x,y,p)&=\sum_{\alpha=0}^{+\infty}J_{\alpha}(x,y,p),
\end{aligned}
\end{equation*}
and similarly for $(f(x),g(x,y))\in\mathcal{F}'$:
\begin{equation*}
\begin{aligned}
f(x) &= \sum_{\alpha=2}^{\infty}f_{\alpha}(x),\\
g(x,y) &= \sum_{\alpha=2}^{\infty}g_{\alpha}(x,y).
\end{aligned}
\end{equation*}
with 
\begin{equation}\label{F'-2}
\begin{aligned}
f_{i}(x) &=f_{i}x^{i},\\
 g_{2}(x,y) &=g_{2,0}x^{2},\\
 g_{3}(x,y) &=g_{3,0}x^{3},\\
 g_{i}(x,y) &=\sum_{\alpha+2\beta=i}g_{\alpha,\beta}x^{\alpha}y^{\beta} \qquad (i\geqslant 4). 
\end{aligned}
\end{equation}
Thus $g_{2,y}\equiv 0$; while $g_{2,x,x}\equiv 0$ implies that $g_{2}(x,y)\equiv 0$. Moreover, $g_{3,y}\equiv 0$; while $g_{3,0}\equiv 0$ implies that $g_{3}(x,y)\equiv 0$.

After these preparations, the fundamental identity \eqref{FE} can be stratified according to various weights as follows. For terms of weight $\alpha$, obviously the starting expansion is:
\[
K_{\alpha}=J_{\alpha}+\cdots. 
\]
For the boxed terms, observe that given a semi-homogeneous function $h_{\alpha}(x,y)$ of order $\alpha$, the term
\[
p^{k}\partial_{x}^{i}\partial_{y}^{j}h_{\alpha}
\]
is semi-homogeneous of order $\alpha-i-2j+k$. Hence, 
\[
K_{\alpha}=J_{\alpha}-\big(g_{\alpha+2,xx}+(2g_{\alpha+2,xy}-f_{\alpha+1,xx})p+g_{\alpha+2,yy}p^{2}\big)+\cdots.
\]
The rest of the terms are non-linear products of $f$, $g$, $J$,  and their derivatives, and they involve only the lower order terms. Thanks to the definition of semi-homogeneity, we are able to stratify our computations in the following way, which will not be possible if the weights are chosen in a different way:
\[
K_{\alpha}=J_{\alpha}-\big(g_{\alpha+2,xx}+(2g_{\alpha+2,xy}-f_{\alpha+1,xx})p+g_{\alpha+2,yy}p^{2}\big)+R_{\alpha}\big(f_{2\leqslant \beta\leqslant \alpha}; g_{2\leqslant \beta\leqslant \alpha+1}; J_{0\leqslant \beta\leqslant \alpha-1}\big).
\]

We write down the first few terms:
\begin{equation}
\begin{aligned}
K_{0} &= J_{0}-g_{2,xx},\\
K_{1} &=J_{1}-\big(g_{3,xx}-f_{2,xx}p\big)+R_{1}(g_{2};J_{0}),\\
K_{2} &=J_{2}-\big(g_{4,xx}+(2g_{4,xy}-f_{3,xx})p+g_{4,yy}p^{2}\big)+R_{2}(f_{2};g_{2\leqslant \beta\leqslant 3};J_{0\leqslant \beta\leqslant 1}),\\
K_{3} &=J_{3}-\big(g_{5,xx}+(2g_{5,xy}-f_{4,xx})p+g_{5,yy}p^{2}\big)+R_{3}(f_{2\leqslant \beta\leqslant 3};g_{2\leqslant \beta\leqslant 4};J_{0\leqslant \beta\leqslant 2}).
\end{aligned}
\end{equation}
We can see that at stage $\alpha$, the terms $f_{\beta\leqslant \alpha}$, $g_{\beta\leqslant \alpha+1}$, $J_{\beta\leqslant \alpha-1}$ have been normalised in the previous $\alpha-1$ steps, and thus will not be free for future normalisations. The only free parameters come from $g_{\alpha+2}$, $f_{\alpha+1}$, which will be used via $L(f_{\alpha+1}, g_{\alpha+2})$  to normalise 
\[J_{\alpha}+R_{\alpha}(f_{2\leqslant \beta\leqslant \alpha};g_{2\leqslant\beta\leqslant \alpha+1};J_{0\leqslant \beta\leqslant \alpha-1}).\]

We will show the following 

\begin{Theorem}
 For $(f,g)\in \mathcal{F}'$, the only fibre-preserving map $x\longmapsto x+f(x)$, $y\longmapsto y+g(x,y)$ that sends the first normal form 
 \[
 y_{xx}=K(x,y,p)
 \]
 to the second normal form 
 \[
 y_{xx}=K'(x,y,p)
 \]
 is the identity map, and consequently $K(x,y,p)=K'(x,y,p)$.
  \end{Theorem}
 
 \begin{proof}
 We will show that such a diffeomorphism implies that $f(x)\equiv 0$ and $g(x,y)\equiv 0$. We compute the first few terms:
 \begin{equation}\label{thm-unique}
 \begin{aligned}
K_{0} &=K_{0}'-g_{2,xx},\\
K_{1} &=K_{1}'-\big(g_{3,xx}-f_{2,xx}p\big)+R_{1}(g_{2};K_{0}'),\\
K_{2} &=K_{2}'-\big(g_{4,xx}+(2g_{4,xy}-f_{3,xx})p+g_{4,yy}p^{2}\big)+R_{2}(f_{2};g_{2},g_{3};K_{0\leqslant \beta\leqslant 1}'),\\
K_{3} &=K_{3}'-\big(g_{5,xx}+(2g_{5,xy}-f_{4,xx})p+g_{5,yy}p^{2}\big)
+R_{3}(f_{2\leqslant \beta\leqslant 3};g_{2\leqslant\beta\leqslant 4};K_{0\leqslant\beta\leqslant 2}'),\\
&\vdots
 \end{aligned}
 \end{equation}
 Using equation \eqref{F'-2}, and the fact that $K_{0}=K_{0}'=0$ for normal forms, the first line in equation \eqref{thm-unique} becomes:
 \[
 0=g_{2,0},
 \]
which implies that $g_{2}(x,y)\equiv 0$. Hence the remainder term in the second line vanishes identically:
 \[
 R_{1}(g_{2};K_{0}')\equiv 0. 
 \]
 Using equation \eqref{F'-2} in the second line, and the fact that $K_{1}\equiv 0$, $K_{1}'\equiv 0$ for normal forms, the equation for the weight 2 component becomes:
\begin{equation*}
\begin{aligned}
0 &= -\big(g_{3,xx}-f_{2,xx}p\big)\\
&= -6g_{3,0}\ x+2f_{2}\ p.
\end{aligned}
\end{equation*}
Consequently, $g_{3}(x,y)\equiv 0$ and $f_{2}(x)\equiv 0$. The remainder term in the third line identically vanishes:
\[
R_{2}(f_{2};g_{2},g_{3};K_{0}',K_{1}')\equiv 0.
\]
Since $K_{2}\equiv 0$ and $K_{2}'\equiv 0$ in normal forms, the weight 2 component of equation \eqref{thm-unique} becomes:
\begin{equation*}
\begin{aligned}
0 &= -\big(4g_{3,xx}+(2g_{4,xy}-f_{3,xx})p+g_{4,yy}p^{2}\big)\\
&= 12g_{4,0}x^{2}+2g_{2,1}y+
(2g_{2,1}x-6f_{3}x)p+2g_{0,2}p^{2},
\end{aligned}
\end{equation*}
which clearly implies that $g_{4}(x,y)\equiv 0$ and $f_{3}(x)\equiv 0$. We will finish the rest of the proof by induction on the weights. For $\alpha\geqslant 3$, assume that at stage $\alpha-1$, we have after normalisation of $K_{\alpha-1}$ and $K_{\alpha-1}'$ the fibre-preserving maps:
\begin{equation}\label{induct}
x\longmapsto x+f_{\alpha+1}+f_{\alpha+2}+\cdots,\qquad
y\longmapsto y+g_{\alpha+2}+g_{\alpha+3}+\cdots.
\end{equation}
The normalisation equation for weight $\alpha$ is:
\begin{equation*}
\begin{aligned}
K_{\alpha} &= K_{\alpha}'-
\big(g_{\alpha+2,xx}+(2g_{\alpha+2,xy}-f_{\alpha+1,xx})p+g_{\alpha+2,yy}p^{2}\big)\\
&\hspace{1.5cm}
+
R_{\alpha}(f_{2},\dots,f_{\alpha};\ g_{2},\dots,g_{\alpha+1};K_{0}',\dots,K_{\alpha-1}').
\end{aligned}
\end{equation*}
We will show that $R_{\alpha}\equiv 0$ by computing the lowest order of the semi-homogenous terms that can occur in $R(f,g,J)$ of the fundamental identity. We substitute the induction hypothesis \eqref{induct} in \eqref{FE}, while writing only terms of lowest order of the components in each summand of $R(f,g,J)$,  we receive:
\begin{equation*}
\begin{aligned}
K(x,y,p) &= 
K'(x,y,p)
- 
\underbrace{L(f_{\alpha+1},g_{\alpha+2})}_{\text{weight }\alpha} -\sum_{\beta\geqslant \alpha}\underbrace{L(f_{\beta+2},g_{\beta+3})}_{\text{weight }\geqslant \alpha+1}\\
&\hspace{0.5cm}
+
\underbrace{\big(
K_{\alpha}'-L(f_{\alpha+1},g_{\alpha+2})\big)}_{\text{weight } \alpha}
\underbrace{\big(-f_{\alpha+1,x}-g_{\alpha+2,y}-g_{\alpha+2,y}f_{\alpha+1,x}\big)}_{\text{weight }\geqslant \alpha}\\
&\hspace{0.5cm}
+
\big(
\underbrace{
-g_{\alpha+2,xx}f_{\alpha+1,x}
-g_{\alpha+2,x}f_{\alpha+1,xx}
-
(2g_{\alpha+2,xy}f_{\alpha+1,x}-g_{\alpha+2,y}f_{\alpha+1,xx}\big)p}_{\text{weight }\geqslant 2\alpha}\\
&\hspace{1.2cm}
\underbrace{-g_{\alpha+2,yy}f_{\alpha+1,x}p^{2}
+
K_{\alpha}'f_{\alpha+1,x}
+K_{\alpha,x}'f_{\alpha+1}+K_{\alpha,x}'g_{\alpha+2}+K_{\alpha,p}'\cdot (P-p)}_{\text{weight }\geqslant 2\alpha}\\
&\hspace{1.2cm}
+\underbrace{(K_{\alpha,x}'f_{\alpha+1}+K_{\alpha,x}'g_{\alpha+2}+K_{\alpha,p}'\cdot (P-p))f_{\alpha+1,x}\big)}_{\text{weight }\geqslant 2\alpha}\\
&\hspace{0.5cm}
+\text{terms of weight }\geqslant 2\alpha,
\end{aligned}
\end{equation*}
while observing that  $P-p$ starts with weight $\alpha+1$. By close inspection, the remainder $R(f,g,J)$ starts with at least order $2\alpha$ which is greater than $\alpha$ since $\alpha\geqslant 3$, and thus $R_{\alpha}\equiv 0$. If we let $\pi_{\alpha}$ be the projection onto the terms of weight $\alpha$, then applying this operator to both sides above, we see that:
\begin{equation*}
\begin{aligned}
K_{\alpha} &= \pi_{\alpha}(K(x,y,p))\\
&=
\pi_{\alpha}\big(K'(x,y,p)-L(f_{\alpha+1},g_{\alpha+2})+...\big)\\
&= K_{\alpha}'-L(f_{\alpha+1},g_{\alpha+2}),
\end{aligned}
\end{equation*}
or,
\[
L(f_{\alpha+1},g_{\alpha+2})
=
K_{\alpha}'-K_{\alpha}\in\mathcal{N}.
\]
The fact $\Im (L)\cap \mathcal{N}=\{0\}$ implies that $L(f_{\alpha+1},g_{\alpha+2})=0$. Since $(f,g)\in\mathcal{F}'$, one also has $(f_{\alpha+1},g_{\alpha+2})\in\mathcal{F}'$. By injective property of $L|_{\mathcal{F}'}$, we have $f_{\alpha+1}\equiv 0$ and $g_{\alpha+2}\equiv 0$.  This finishes the induction process, and the proof is completed. 
 \end{proof}



\section{Proof of the Existence of Normal Forms}

\subsection{} The existence problem asks if there is a convergent (analytic) fibre-preserving map $\Phi$:
\[
X=\varphi(x),\qquad Y=\psi(x,y),
\]
sending $Y_{XX}=J(X,Y,P)$ to $y_{xx}=K(x,y,p)$ which is in normal form. In practice, this is delicate to answer because the maps can be very non-linear in nature. A good strategy is to search for $\Phi$ as a composition of four simpler fibre-preserving maps:
\[
\Phi=\Phi_{1}\circ \Phi_{2}\circ \Phi_{3}\circ \Phi_{4},
\]
where the first transformation $\Phi_{1}$ sends $Y_{XX}=J(X,Y,P)$ to $y_{xx}=J_{1}(x,y,p)$ with the property D1, which is $J_{1}(x,y,0)\equiv 0$. Next the second map $\Phi_{2}$ sends $Y_{XX}=J_{1}(X,Y,P)$ to $y_{xx}=J_{2}(x,y,p)$ with the property D2 while still keeping the property D1, meaning $J_{2}(x,y,0)\equiv 0$ and $J_{2,p}(0,y,0)\equiv 0$. The third map $\Phi_{3}$ sends $Y_{XX}=J_{2}(X,Y,P)$ to $y_{xx}=J_{3}(x,y,p)$ with properties D1, D2, D3, that is $J_{3}(x,y,0)\equiv 0$, $J_{3,p}(0,y,0)\equiv 0$, and $J_{3,p}(x,0,0)\equiv 0$. Finally $\Phi_{4}$ sends $Y_{XX}=J_{3}(X,Y,P)$ to the normal form. 

A pertinent question therefore is: \textit{how can we find such $\Phi_{i}$}? Fortunately we can make a good guess based on the formal computations in the preceding section. For example, to find $\Phi_{2}$ giving us D2, the results from  formal computations in Theorem \ref{thm-nf} give us an explicit fixed pair of solutions $f(x)$, $g(x,y)$. The flow at $t=1$ of the vector field $X=f(x)\partial_{x}+g(x,y)\partial_{y}$ (i.e. $\exp(X)$), which is 
\begin{equation}
\begin{aligned}
X &= x+f(x)+\cdots,\\
Y &= y+g(x,y)+\cdots,
\end{aligned}
\end{equation}
should in principle be our desired map $\Phi_{2}$. However in reality, such a  transformation, while achieving D2, will break D1. Thus our guess is based on modifying $\exp(X)$ via
\begin{equation}
\begin{aligned}
X &= x+f(x)+a(x),\\
Y &= y+g(x,y)+b(x,y).
\end{aligned}
\end{equation}
The partial differential equations corresponding to D1, D2, will provide a system of differential equations in $f(x)$, $a(x)$, $g(x,y)$, $b(x,y)$. The existence and uniqueness of the solutions are then guaranteed by the Cauchy-Kovalevskaya theorem, since we are working in the real-analytic category. This idea will be the main guiding principle for the proof of the existence of normal forms. 

\subsection{The condition D1.} After explaining our strategy, we proceed to find $\Phi_{1}$ to satisfy D1, which is equivalent to $K_{i,0}(y)=0$ for all $i\geqslant 0$. We look at equation \eqref{sol-1} in Theorem \ref{thm-nf}, which tells us that we should be looking at the following infinitesimal symmetry:
\begin{equation}\label{X-D1}
\begin{aligned}
{\sf X} &= \bigg(\sum_{i=0}^{\infty}g_{i+2}(y)\bigg)\partial_{y}\\
&=
\bigg(
\sum_{i=0}^{\infty}\frac{-1}{(i+1)(i+2)}J_{i,0}(y)x^{i+2}\bigg)\partial_{y}.
\end{aligned}
\end{equation}
To simplify our discussion, we rewrite this vector field as
\[
{\sf X}=g(x,y)\partial_{y},
\]
and since the coefficient in front of $\partial_{y}$ in \eqref{X-D1} starts with $x^{2}$, we require that 
\[
g(0,y)\equiv 0, \qquad g_{x}(0,y)\equiv 0.
\]
Computing the flow  $\exp(X)$ of this vector field at $t=1$, the fibre-preserving maps become 
\begin{equation}
\begin{aligned}
X &= x,\\
Y &= y+G(x,y).
\end{aligned}
\end{equation}
Moreover, we claim that $G(0,y)\equiv 0$ and $G_{x}(0,y)\equiv 0$. This is because looking at the equation of the flow:
\begin{equation}
\begin{aligned}
\gamma_{1}'(t) &= 0,\\
\gamma_{2}'(t) &= g(\gamma_{1}(t),\gamma_{2}(t)),
\end{aligned}
\end{equation}
along with the initial conditions $\gamma_{1}(0)=x$, $\gamma_{2}(0)=y$, we have the convergent expansion:
\begin{equation}
\begin{aligned}
\gamma_{2}(t) &= 
\gamma_{2}(0)+\gamma_{2}'(0)t + \gamma_{2}^{(2)}(0)\tfrac{t^{2}}{2!}+\gamma_{2}^{(3)}(0)\tfrac{t^{3}}{3!}+\cdots\\
&= y+g(x,y)t+\gamma_{2}^{(2)}(0)\tfrac{t^{2}}{2!}+\gamma_{2}^{(3)}(0)\tfrac{t^{3}}{3!}+\cdots.
\end{aligned}
\end{equation}
It remains to show that $\gamma_{2}^{(k)}(t)=:H_{k}(\gamma_{1}(t),\gamma_{2}(t))$ has the same property of being $O(x^2)$ at $t=0$.  We will prove this by induction by supposing that this is true at stage $k$. Then 
\begin{equation}
\begin{aligned}
\gamma_{2}^{(k+1)}(t) &=H_{k,x}(\gamma_{1}(t),\gamma_{2}(t))\gamma_{1}'(t)+H_{k,y}(\gamma_{1}(t),\gamma_{2}(t))\gamma_{2}'(t)\\
&= H_{k,y}(\gamma_{1}(t),\gamma_{2}(t))g(\gamma_{1}(t),\gamma_{2}(t)).
\end{aligned}
\end{equation}
At $t=0$, we get $\gamma_{2}^{(k+1)}(0)=H_{k,y}(x,y)g(x,y)$. Since $g(x,y)=O(x^2)$, so is $\gamma_{2}^{(k+1)}(0)$, and we have proved our claim.  These two properties of $G(x,y)$ will provide initial conditions for some Cauchy-Kovalevskaya system. 
\begin{Proposition}
There exists a unique map
\begin{equation}
\begin{aligned}
X &= x,\\
Y &= y+G(x,y),
\end{aligned}
\end{equation}
with $G(0,y)\equiv 0$ and $G_{x}(0,y)\equiv 0$, sending $Y_{XX}=J(X,Y,P)$ to $y_{xx}=K(x,y,p)$ with $K(x,y,0)\equiv 0$.
\end{Proposition}
\begin{proof}
Based on the formal computations, we would expect that the D1 condition requires a second order partial differential equation in $x$ of $G(x,y)$. Indeed, we prolong this map to second order,
\begin{equation}
\begin{aligned}
P &= \frac{{\sf D}_{x}(Y)}{{\sf D}_{x}(X)}=p(1+G_{y})+G_{x},\\
Y_{XX} &= \frac{{\sf D}_{x}P}{{\sf D}_{x}X}=y_{xx}(1+G_{y})+G_{xx}+2pG_{xy}+p^{2}G_{yy}.
\end{aligned}
\end{equation}
Substituting these into $Y_{XX}=J(X,Y,P)$ and solving for $y_{xx}$, the resulting second order ordinary differential equation $y_{xx}=K(x,y,p)$ becomes
\begin{equation}
\begin{aligned}
y_{xx} &= \frac{1}{1+G_{y}}\big[
J\big(x,G(x,y),p(1+G_{y}(x,y))+G_{x}(x,y)\big)-G_{xx}(x,y)-2pG_{xy}(x,y)-p^{2}G_{yy}(x,y)\big]\\
&=:K(x,y,p).
\end{aligned}
\end{equation}
Then the requirement $K(x,y,0)\equiv 0$ holds if and only if we have the following Cauchy-Kovalevskaya system:
\[
J\big(x,G(x,y),G_{x}(x,y)\big)-G_{xx}(x,y)\equiv 0,
\]
along with the assumptions $G(0,y)\equiv 0$ and $G_{x}(0,y)\equiv 0$. The existence and uniqueness of a solution follows from the Cauchy-Kovalevskaya theorem.
\end{proof}

\subsection{The conditions D1 and D2.} We will now assume that $Y_{XX}=J(X,Y,P)$ satisfies the D1 condition $J(X,Y,0)\equiv 0$. We seek a second fibre-preserving map $\Phi_{2}$ that brings it to $y_{xx}=K(x,y,p)$ satisfying both the D1 and D2 conditions, which are $K(x,y,0)\equiv 0$ and $K_{p}(0,y,0)\equiv 0$. Recall that the D2 condition is equivalent to $K_{0,1}(y)\equiv 0$. Based on equations \eqref{sol-2-f}, \eqref{sol-2-g} in Theorem \ref{thm-nf}, we should be looking at the flow of the following infinitesimal symmetry
\[
{\sf X}=\big(f_{2}x^{2}\big)\partial_{x}+\big(g_{1}(y)x\big)\partial_{y},
\]
where 
\begin{equation}\label{D2-eq1}
f_{2}=-\tfrac{1}{4}J_{0,1}(0),\qquad 
g_{1}'(y)=\tfrac{1}{2}\big(-J_{0,1}(y)+J_{0,1}(0)\big).
\end{equation}
It is therefore clear that $g_{1}'(0)=0$. Moreover, by assumption in the introduction, $g_{1}(0)=0$. The next proposition shows the nature of $\exp({\sf X})$. 
\begin{Lemma}\label{lem-D2}
For any $s(y)$, the flow of ${\sf X}=f_{2}x^{2}\partial_{x}+s(y)x\partial_{y}$ at $t=1$, $\exp({\sf X})$, defines the fibre preserving map of the form:
\begin{equation}
\begin{aligned}
X &= \frac{x}{1-f_{2}x},\\
Y &= y+s(y)x+h(x,y),
\end{aligned}
\end{equation}
for some unique $h(x,y)$ with $h(0,y)\equiv 0$ and $h_{x}(0,y)\equiv 0$.
\end{Lemma}
\begin{proof}
The flow $\exp(tX)$ is described by the system of ordinary differential equations:
\begin{equation*}
\begin{aligned}
\gamma_{1}'(t) &= f_{2}\gamma_{1}(t)^{2},\\
\gamma_{2}'(t) &= s(\gamma_{2}(t))\gamma_{1}(t),
\end{aligned}
\end{equation*}
subject to the initial conditions $\gamma_{1}(0)=x$, $\gamma_{2}(0)=y$. Solving the first equation is straightforward:
\[
\gamma_{1}(t)=\frac{x}{1-f_{2}xt}=O(x).
\]
Solving the second order ordinary equation with a closed formula is probably not possible, but at least there exists a convergent power series solution in $t$:
\begin{equation}
\begin{aligned}
\gamma_{2}(t) 
&= 
y + s(y)x t + \gamma_{2}''(0) \tfrac{t^{2}}{2!}+
\gamma_{2}'''(0)\tfrac{t^{3}}{3!}
+
\cdots.
\end{aligned}
\end{equation}
We claim that for any $k\geqslant 1$,
\[
\frac{d^{k}}{dt^{k}}\gamma_{2}'(t)=O(\gamma_{1}(t)^{k+1}).
\]
Indeed if at stage $k$, we have for some function $H(\lambda,\mu)$ such that 
\[
\frac{d^{k}}{dt^{k}}\gamma_{2}'(t)=H_{k}(\gamma_{1}(t),\gamma_{2}(t))\gamma_{1}(t)^{k+1},
\]
then differentiating both sides with respect to $t$, and using the ordinary differential equations to make appropriate replacements of $\gamma_{1}'$ and $\gamma_{2}'$:
\begin{equation}
\begin{aligned}
\frac{d^{k+1}}{dt^{k+1}}\gamma_{2}'(t)
&=
H_{k,\lambda}\big(\gamma_{1}(t),\gamma_{2}(t)\big)\gamma_{1}'(t)\gamma_{1}(t)^{k+1}
+
H_{k,\mu}\big(\gamma_{1}(t),\gamma_{2}(t)\big)\gamma_{2}'(t)\gamma_{1}(t)^{k+1}\\
&\hspace{0.5cm}
+
(k+1)H_{k}\big(\gamma_{1}(t),\gamma_{2}(t)\big)\gamma_{1}(t)^{k}\gamma_{1}'(t)\\
&=
H_{k,\lambda}\big(\gamma_{1}(t),\gamma_{2}(t)\big)f_{2}\gamma_{1}(t)^{k+3}
+
H_{k,\mu}\big(\gamma_{1}(t),\gamma_{2}(t)\big)s(\gamma_{2}(t))\gamma_{1}(t)^{k+2}\\
&\hspace{0.5cm}
+
(k+1)H_{k}\big(\gamma_{1}(t),\gamma_{2}(t)\big)f_{2}\gamma_{1}(t)^{k+2}\\
&=:
H_{k+1}\big(\gamma_{1}(t),\gamma_{2}(t)\big)\gamma_{1}(t)^{k+2},
\end{aligned}
\end{equation}
which finishes the induction. At $t=0$, the initial conditions give
\begin{equation*}
\begin{aligned}
\gamma_{2}^{(k+2)}(0)
=
H_{k+1}\big(x,y\big)x^{k+2}
=
O(x^{k+2}).
\end{aligned}
\end{equation*}
Hence the diffeomorphism at $t=1$ becomes
\begin{equation}
\begin{aligned}
X &= \frac{x}{1-f_{2}x},\\
Y &=y + s(y)x  + \gamma_{2}''(0) \tfrac{1}{2!}+
\gamma_{2}'''(0)\tfrac{1}{3!}
+
\cdots\\
&= 
y+s(y)x+h(x,y),
\end{aligned}
\end{equation}
where $h(0,y)\equiv 0$ and $h_{x}(0,y)\equiv 0$. This finishes the proof. 
\end{proof}

By guessing based on Lemma \ref{lem-D2}, our next proposed fibre-preserving maps to obtain D2 while still preserving D1 will therefore be of the form
\begin{equation}
\begin{aligned}
X &= \frac{x}{1-f_{2}x},\\
Y &= y + s(y)x + h(x,y),
\end{aligned}
\end{equation}
where $s(0)=0$, $h(0,y)\equiv 0$, and  $h_{x}(0,y)\equiv 0$. The quantities $f_{2}$, $s(y)$, $h(x,y)$ will be determined later based on the needs of our normalisation. To simplify matters, we would like to have $f_{2}=0$. From the formal computations, this holds only after we have $J_{0,1}(0)=0$, or equivalently  $J_{p}(0,0,0)=0$ (see equation \eqref{D2-eq1}). Thus we start by:
\begin{Lemma}
Let $Y_{XX}=J(X,Y,P)$ be a second order ordinary differential equation with $J(X,Y,0)\equiv 0$. Then for any $f_{2}$, the map 
\begin{equation}
\begin{aligned}
X &=\frac{x}{1-f_{2}x},\\
Y &= y,
\end{aligned}
\end{equation}
sends it to $y_{xx}=K(x,y,p)$ with the same D1 condition $K(x,y,0)\equiv 0$. Moreover, there exists a unique $f_{2}$ such that $K_{p}(0,0,0)=0$. 
\end{Lemma}
\begin{proof}
We prolong the map to second order
\begin{equation}
\begin{aligned}
P &= (1-f_{2}x)^{2}p,\\
Y_{XX} &= (1-f_{2}x)^{3}\big(y_{xx}(1-f_{2}x)-2pf_{2}\big),
\end{aligned}
\end{equation}
and by substitution into $Y_{XX}=J(X,Y,P)$, we obtain the new second order ordinary differential equation:
\begin{equation*}
\begin{aligned}
y_{xx} &= 
\frac{1}{(1-xf_{2})^{4}}
\bigg[
\big(-2f_{2}^{4}x^{3}+6f_{2}^{3}x^{2}-6f_{2}^{2}x+2f_{2}\big)p
+
J\bigg(\frac{x}{1-f_{2}x},y,p(1-f_{2}x)^{2}\bigg)\bigg]\\
&=:K(x,y,p).
\end{aligned}
\end{equation*}
Hence at $p=0$, 
\[
K(x,y,0)=\frac{1}{1-xf_{2}}J\bigg(\frac{x}{1-f_{2}x},y,0\bigg)\equiv 0,
\]
since $J(x,y,0)\equiv 0$. Thus the condition D1 is preserved. 

Next, a differentiation of $K$ with respect to $p$ yields:
\[
K_{p}(0,0,0)=2f_{2}+J_{p}(0,0,0).
\]
We may choose $f_{2}:=-\tfrac{1}{2}J_{p}(0,0,0)$ to normalise $K_{p}(0,0,0)$ to zero, and the proof is completed. 
\end{proof}

Henceforth we will assume that $J_{p}(0,0,0)=0$ so that $f_{2}$ will no longer be needed. We are now in a position to complete the second normalisation. 

\begin{Proposition}
Let $Y_{XX}=J(X,Y,P)$ be a second order ordinary differential equation, with $J(X,Y,0)\equiv 0$ and $J_{P}(0,0,0)\equiv 0$. Then there exists a unique fibre preserving transformation of the form 
\begin{equation}
\begin{aligned}
X &= x,\\
Y &= y + s(y)x + h(x,y),
\end{aligned}
\end{equation}
with $h(0,y)\equiv 0$, $h_{x}(0,y)\equiv 0$ (imitating $\exp({\sf X})$), and $s(0)=0$ (from $g_{1}(0)=0$), sending it to $y_{xx}=K(x,y,p)$ satisfying both D1 and D2 conditions. In other words, $K(x,y,0)\equiv 0$ and $K_{p}(0,y,0)\equiv 0$.
\end{Proposition}
\begin{proof}
As before, we prolong the map to second order:
\begin{equation}
\begin{aligned}
P &= (1+h_{y}+s_{y}(x))p+s+h_{x},\\
Y_{XX} &=y_{xx}(1+h_{y}+xs_{y})+h_{xx}+2h_{xy}p+(h_{yy}+s_{y}x)p^{2},
\end{aligned}
\end{equation}
and by substituting in $Y_{XX}=J(X,Y,P)$, we obtain
\begin{equation*}
\begin{aligned}
K(x,y,0) &= \frac{1}{1+s_{y}x+h_{y}}
\bigg[
J\big(x,y+s(y)x+h(x,y),s(y)+h_{x}(x,y)\big)-h_{xx}(x,y)\bigg],\\
K_{p}(0,y,0) &= -2s_{y}(y)+J_{P}(0,y,s(y)).
\end{aligned}
\end{equation*}
To complete the normalisation, we need both of them to identically vanish. Using the initial condition $s(0)=0$, we may solve the first order ordinary differential equation in the second line to obtain the solution for $s(y)$, which will be substituted into the first equation for $K(x,y,0)$. The resulting first equation, being a second order partial differential equation in $x$ of $h(x,y)$, can be solved with the initial conditions $h(0,y)\equiv 0$, $h_{x}(0,y)\equiv 0$ thanks to the Cauchy-Kovalevskaya theorem. The proof is completed. 
\end{proof}

\subsection{The conditions D1, D2, and D3.} After completing our first and second normalisations, we are given $Y_{XX}=J(X,Y,P)$ with $J(X,Y,0)\equiv 0$ and $J_{P}(0,Y,0)\equiv 0$. We are going to find $\Phi_{3}$ bringing it to $y_{xx}=K(x,y,p)$ satisfying D1, D2, D3. The condition D3 is equivalent to $K_{i,1}(0)=0$ for all $i\geqslant 1$. Based on equation \eqref{sol-3} in Theorem \ref{thm-nf}, we should be looking at the following fibre-preserving symmetry:
\[
{\sf X}=
\bigg(\sum_{j=3}^{\infty}f_{j}x^{j}\bigg)\partial_{x},
\]
where for $i\geqslant 1$:
\[
f_{i+2}=\frac{1}{(i+1)(i+2)}\big(J_{i,1}(0)-\tfrac{2}{i}J_{i-1,0}'(0)\big).
\]
The fibre-preserving map $\exp({\sf X})$ associated to the vector field is of the form:
\begin{equation*}
\begin{aligned}
X &= x+\bigg(\sum_{j=3}^{\infty}f_{j}x^{j}\bigg)+O(x^{5}),\\
Y &= y. 
\end{aligned}
\end{equation*}
As previously mentioned, this choice may break either D1 or D2 condition, which necessitates the modification of $X$ with a function $a(x)$ of order $x^{5}$ to preserve them:
\begin{equation*}
\begin{aligned}
X &= x+\bigg(\sum_{j=3}^{\infty}f_{j}x^{j}\bigg)+O(x^{5}) + a(x),\\
Y &= y. 
\end{aligned}
\end{equation*}
As a result,  our guess of our next normalisation will be of the form:
\begin{equation}
\begin{aligned}
X &= x+f(x),\\
Y &= y,
\end{aligned}
\end{equation}
with $f(0)=0$, $f'(0)=0$, $f''(0)=0$, since the coefficient of $\partial_{x}$ in ${\sf X}$ starts with $x^{3}$. We will be expecting a third order ordinary differential equation for $f(x)$. 
\begin{Proposition}
Let $Y_{XX}=J(X,Y,P)$ be a second order ordinary differential equation satisfying  conditions D1 and D2, that is $J(X,Y,0)\equiv 0$ and $J_{P}(0,Y,0)\equiv 0$. Then any fibre preserving mapping of the form
\begin{equation}
\begin{aligned}
X &= x+f(x),\\
Y &= y,
\end{aligned}
\end{equation}
with $f(0)=0$, $f'(0)=0$, $f''(0)=0$, sends it to $y_{xx}=K(x,y,p)$ with the same conditions D1, D2, being satisfied, that is $K(x,y,0)\equiv 0$ and $K_{p}(0,y,0)\equiv 0$. Moreover, there exists a unique $f(x)$ such that $K_{p}(x,0,0)\equiv 0$, satisfying D3. 
\end{Proposition}
\begin{proof}
We  first show that any such map with $f(0)=f'(0)=f''(0)=0$ preserves D1 and D2. Prolonging the map to second order,
\begin{equation}
\begin{aligned}
P &= \frac{p}{1+f_{x}},\\
Y_{XX} &= \frac{y_{xx}(1+f_{x})-pf_{xx}}{(1+f_{x})^{3}},
\end{aligned}
\end{equation}
and substituting these into $Y_{XX}=J(X,Y,P)$, we obtain $y_{xx}=K(x,y,p)$ where 
\begin{equation}
\begin{aligned}
K(x,y,p)&=\frac{1}{1+f_{x}}
\bigg[
\big(f_{x}^{3}+3f_{x}^{2}+3f_{x}+1\big)J\bigg(x+f,y,\frac{p}{1+f_{x}}\bigg)
+
pf_{xx}
\bigg].
\end{aligned}
\end{equation}
Then at $p=0$, 
\begin{equation}
\begin{aligned}
K(x,y,0)&=\frac{1}{1+f_{x}}
\bigg[
\big(f_{x}^{3}+3f_{x}^{2}+3f_{x}+1\big)J\big(x+f,y,0\big)\bigg]\\
&\equiv 0,
\end{aligned}
\end{equation}
 vanishes since $J(x,y,0)\equiv 0$. 
 
 Next, we verify that the condition D2 is preserved. We differentiate $K(x,y,p)$ with respect to $p$:
\begin{equation*}
\begin{aligned}
K_{p}(x,y,p)&=\frac{1}{(1+f_{x})^{2}}
\bigg[
\big(f_{x}^{3}+3f_{x}^{2}+3f_{x}+1\big)J_{p}\bigg(x+f,y,\frac{p}{1+f_{x}}\bigg)
+
pf_{xx}(1+f_{x})
\bigg].
\end{aligned}
\end{equation*}
Using $f(0)=0$, $f_{xx}(0)=0$ and $J_{p}(0,y,0)\equiv 0$, we see that $K_{p}(0,y,0)\equiv 0$, and thus D2 is again satisfied. In conclusion, both D1 and D2 are preserved. 

Now to obtain D3 for $K(x,y,p)$, we need to find a third order ordinary differential equation for $f(x)$ so that $K_{p}(x,0,0)\equiv 0$. From the preceding equation, we see that $K_{p}(x,y,x)$ contains $f$ up to its second order derivative $f_{xx}$, and hence we differentiate it with respect to $x$ and look at $K_{px}(x,0,0)\equiv 0$. This is sufficient to obtain $K_{p}(x,0,0)\equiv 0$ since we have both $K_{p}(0,0,0)\equiv 0$ and $K_{px}(x,0,0)\equiv 0$:
\begin{equation*}
\begin{aligned}
K_{px}(x,0,0)
&=
\frac{f_{xxx}}{1+f_{x}}
-
\frac{f_{xx}^{2}}{(1+f_{x})^{2}}
+
(1+f_{x})^{2}J_{xp}(x+f,0,0)\\
&\hspace{0.5cm}
+f_{xx}J_{p}(x+f,0,0).
\end{aligned}
\end{equation*}
For this to identically vanish, it suffices to solve for $f$:
\[
f_{xxx}
=
\frac{f_{xx}^{2}}{1+f_{x}}
-
(1+f_{x})^{3}J_{xp}(x+f,0,0)
-
f_{xx}(1+f_{x})J_{p}(x+f,0,0).
\]
Along with the initial conditions $f(0)=0$, $f_{x}(0)=0$, $f_{xx}(0)=0$, this equation has a unique solution for $f$, and the proof is completed. 
\end{proof}

We have therefore completed the third normalisation.

\subsection{The final step: normal form.}

After the third normalisation, we receive $Y_{XX}=J(X,Y,P)$ with $J(X,Y,0)\equiv 0$, $J_{P}(0,Y,0)\equiv 0$ and $J_{P}(X,0,0)\equiv 0$. We are left with the final task of finding $\Phi_{4}$ bringing it to the normal form $y_{xx}=K(x,y,p)$ satisfying D1, D2, D3, D4. The D4 condition is equivalent to $K_{0,2}(y)\equiv 0$. Based on equation \eqref{sol-4} in Theorem \ref{thm-nf}, we should be looking at the following fibre-preserving infinitesimal symmetry:
\[
{\sf X}=g_{0}(y)\partial_{y},
\]
where $g_{0}(0)=0$, $g_{0}'(0)=0$, and $g_{0}''(y)=-J_{0,2}(y)$. Thus we have a good guess for our final normalisation.
\begin{Proposition}
Let $Y_{XX}=J(X,Y,P)$ be a second order ordinary differential equation with $J(X,Y,0)\equiv 0$, $J_{P}(0,Y,0)\equiv 0$, and $J_{P}(X,0,0)\equiv 0$. Then any  fibre-preserving transformation of the form
\begin{equation}
\begin{aligned}
 X &= x,\\
 Y &= y+g(y),
\end{aligned}
\end{equation}
with $g(0)=0$, $g_{y}(0)=0$, sends it to $y_{xx}=K(x,y,p)$ with the same properties $K(x,y,0)\equiv 0$, $K_{p}(0,y,0)\equiv 0$, and $K_{p}(x,0,0,)\equiv 0$. Moreover, there exists a unique $g(y)$ such that $K_{pp}(0,y,0)\equiv 0$.
\end{Proposition}
\begin{proof}
Firstly, we show that such a transformation preserves D1, D2, D3. We prolong the map to second order:
\begin{equation}
\begin{aligned}
P &= (1+g_{y})p,\\
Y_{XX} &= y_{xx}(1+g_{y})+g_{yy}p^{2}.
\end{aligned}
\end{equation}
We substitute these into $Y_{XX}=J(X,Y,P)$, and we obtain the second order ordinary differential equation $y_{xx}=K(x,y,p)$ with 
\[
K(x,y,p) = 
\frac{-1}{1+g_{y}}
\bigg[g_{yy}p^{2}-J\big(x,y+g(y),(1+g_{y})p\big)\bigg].
\]
Thus property D1 is preserved:
\[
K(x,y,0)=\frac{-1}{1+g_{y}}\bigg[-J\big(x,y+g(y),0\big)\bigg]\equiv 0.
\]
Next, we differentiate $K$ with respect to $p$:
\[
K_{p}(x,y,p)=
\frac{-1}{1+g_{y}}
\bigg[
2g_{yy}p-J_{p}\big(x,y+g(y),(1+g_{y})p\big)(1+g_{y})\bigg].
\]
Using $g(0)=0$, $g_{y}(0)=0$, it follows that D2 holds:
\[
K_{p}(x,0,0)=-J_{p}(x,0,0)\equiv 0,
\]
and D3 as well:
\[
K_{p}(0,y,0)=J_{p}\big(0,y+g(y),0\big)\equiv 0.
\]

To achieve the final normalisation, we differentiate $K_{p}$ with respect to $p$, we set $x=p=0$, and we obtain
\begin{equation}
\begin{aligned}
K_{p}(0,y,0)=
\frac{-1}{1+g_{y}}
\bigg[2g_{yy}-J_{pp}\big(0,y+g(y),0\big)(1+g_{y})^{2}\bigg].
\end{aligned}
\end{equation}
For $K_{pp}(0,y,0)\equiv 0$ to hold, it suffices to solve for $g$ the following second order ordinary differential equation:
\[
g_{yy}=\tfrac{1}{2}J_{pp}\big(0,y+g(y),0\big)(1+g_{y})^{2}.
\]
Along with the initial conditions $g(0)=0$, $g_{y}(0)=0$, this equation can be solved for $g$ and the solution is unique once again thanks to the Cauchy-Kovalevskaya theorem. The proof is completed. 
\end{proof}

\subsection{The Existence Theorem.} \label{ET}After having completed the four normalisations, we summarise these steps in the following
\begin{Summary}\label{SE}
Let $Y_{XX}=J(X,Y,P)$ be a second order ordinary differential equation. Then there exists a finite sequence of fibre-preserving transformation sending it to a normal form $y_{xx}=K(x,y,p)$. The four steps are as follows:
\begin{enumerate}
\item First, the transformation 
\begin{equation}
\begin{aligned} 
X &= x,\\
Y &=y+G(x,y),
\end{aligned}
\end{equation}
with $G(0,y)\equiv 0$, $G_{x}(0,y)\equiv 0$, satisfying 
\[
G_{xx}(x,y)=J\big(x, G(x,y), G_{x}(x,y)\big),
\]
brings it to $y_{xx}=J_{1}(x,y,p)$ with $J_{1}(x,y,0)\equiv 0$. 

\item Next, given $Y_{XX}=J_{1}(X,Y,P)$ as before, we may perform the transformation
\begin{equation}
\begin{aligned}
X &= \frac{x}{1+\frac{1}{2}J_{1,p}(0,0,0)x},\\
Y &= y,
\end{aligned}
\end{equation}
to bring $J_{1,p}(0,0,0)$ to zero. Once this is done, apply the next fibre-preserving map: 
\begin{equation}
\begin{aligned}
X &= x,\\
Y &= y+s(y)x+h(x,y),
\end{aligned}
\end{equation}
where $s(0)=0$, $h(0,y)\equiv 0$, $h_{x}(0,y)\equiv 0$, satisfying 
\begin{equation}
\begin{aligned}
h_{xx} &= J_{1}(x,y+s(y)x,s+h_{x}),\\
s_{y} &= \frac{1}{2}J_{1,p}(0,y,s(y)),
\end{aligned}
\end{equation}
to bring it to $y_{xx}=J_{2}(x,y,p)$ with $J_{2}(x,y,0)\equiv 0$, and $J_{2,p}(0,y,0)\equiv 0$. 

\item Then given $Y_{XX}=J_{2}(X,Y,P)$, we let 
\begin{equation}
\begin{aligned}
X &= x+f(x),\\
Y &= y,
\end{aligned}
\end{equation}
with $f(0)=0$, $f_{x}(0)=0$, $f_{xx}(0)=0$, satisfying the third order ordinary differential equation
\[
f_{xxx}
=
\frac{f_{xx}^{2}}{1+f_{x}}
-
(1+f_{x})^{3}J_{xp}(x+f,0,0)
-
f_{xx}(1+f_{x})J_{p}(x+f,0,0).
\]

This map sends it to $y_{xx}=J_{3}(x,y,p)$ with $J_{3}(x,y,0)\equiv 0$, $J_{3,p}(0,y,0)\equiv 0$ and $J_{3,p}(x,0,0)\equiv 0$, satisfying D1, D2, D3. 

\item Finally, given $Y_{XX}=J_{3}(X,Y,P)$, we perform the following transformation
\begin{equation}
\begin{aligned}
X &= x,\\
Y &= y+g(y),
\end{aligned}
\end{equation}
with $g(0)=0$, $g_{y}(0)=0$, and 
\[
g_{yy}=\frac{1}{2}J_{3,pp}(0,y+g(y),0)(1+g_{y})^{2},
\]
sending it to the normal form $y_{xx}=K(x,y,p)$. \qed
\end{enumerate}
\end{Summary}

We have therefore proved the following
\begin{Theorem}\label{TE}
There exists a convergent (analytic) fibre-preserving point transformation $x\longmapsto f(x)$, $y\longmapsto g(x,y)$ sending a second order ordinary differential equation $y_{xx}=J(x,y,y_{x})$ to the normal form. 
\end{Theorem}

As a result, we see that if $y_{xx}=K(x,y,p)$ is in normal form, then 
\begin{equation}
\begin{aligned}
I_{1}(K)(0) &=K_{ppp}(0)=K_{0,0,3},\\
I_{2}(K)(0) &= {\sf D}_{x}K_{pp}-K_{py}\big|_{(x,y,p)=0}\\
&=
K_{xpp}(0)+
\underbrace{
p\cdot K_{ypp}+
KK_{pp}
-K_{py}}_{\text{these terms vanish at the origin since }K\text{ is normal}}\big|_{(x,y,p)=0}\\
&=K_{1,0,2},\\
I_{3}(K)(0) &=\underbrace{ {\sf D}_{x}K_{py}}_{=K_{1,1,1}\text{ at the origin}}+\underbrace{K_{pp}K_{y}-K_{py}K_{p}-2K_{yy}}_{\text{these terms vanish at the origin since }K\text{ is normal}}\big|_{(x,y,p)=0}\\
&=K_{1,1,1}.
\end{aligned}
\end{equation}

This may seem like a  coincidence, but in our future manuscript, we will show that the higher order coefficients of a normal form are related to the values of higher order differential invariants at the origin, thanks to Olver's recurrence formulae.

An important application of the existence of the normal form is the following:

\begin{Corollary}
A second order ordinary differential equation admits a six-dimensional Lie group of fibre-preserving point symmetries if and only if there exists a fibre-preserving map that sends it to a normal form:
\[
y'' = M(x,y)p^2+N(x,y)p,
\]
where the functions $M$, $N$ are given by:
\begin{equation}
\begin{aligned}
M(x,y) &= \bigg(J_{1,0,2}x
+\sum_{\substack{i,j\geqslant 0\\ i+j\geqslant 1}}J_{1+i,j,2}x^{1+i}y^{j}\bigg),\\
N(x,y) &= \bigg(J_{1,1,1}xy 
+ \sum_{\substack{i,j\geqslant 0\\ i+j\geqslant 1}}J_{1+i,1+j,1}x^{1+i}y^{1+j}\bigg).
\end{aligned}
\end{equation}
Moreover, the following relations hold:
\begin{equation}
\begin{aligned}
2M_{x} &= N_{y},\\
N_{xy} &= N_{y}N.
\end{aligned}
\end{equation}
Any normal second order ODE that satisfies the criteria above is equivalent to $y''=0$.
\end{Corollary}

\begin{proof}
The existence theorem implies that we may assume $y_{xx}=K(x,y,p)$ is in normal form. Assuming that $I_{1}(K)\equiv 0$, $I_{2}(K)\equiv 0$ and $I_{3}(K)\equiv 0$. The first vanishing condition $I_{1}(K)=K_{ppp}\equiv 0$ implies the existence of functions $M(x,y)$, $N(x,y)$, $C(x,y)$ such that 
\[
y_{xx}=M(x,y)p^{2}+N(x,y)p+C(x,y).
\]
But $C(x,y)\equiv 0$ since $K$ is normal, thus 
\[
y_{xx}=M(x,y)p^{2}+N(x,y)p.
\]
The remaining vanishing conditions $I_{2}(K)\equiv I_{3}(K)\equiv 0$ gives the stated differential relations between $M(x,y)$ and $N(x,y)$. It remains to show that these equations imply that the ODE is equivalent to $y''=0$. 

It suffices to show that in the normal form, the second equation $N_{xy}=N_{y}N$ implies that $N\equiv 0$. To this end, we write $N(x,y)$ as 
\[
N(x,y)=\sum_{k=0}^{\infty}N_{k}(y)x^{k}.
\]
Substituting this into the equation $N_{xy}=N_{y}N$ , we receive
\begin{equation}
\begin{aligned}
\sum_{k=0}^{\infty}N_{k+1}'(y)\cdot (k+1)x^{k}
&=
\bigg(\sum_{k=0}^{\infty}N_{k}'(y)x^{k}\bigg)\bigg(\sum_{k=0}^{\infty}N_{k}(y)x^{k}\bigg)\\
&=
\sum_{k=0}^{\infty}\bigg(\sum_{i+j=k}N_{i}'(y)N_{j}(y)\bigg)x^{k}.
\end{aligned}
\end{equation}
By inspection, immediately we see that for $k=0$:
\[
N_{1}'(y)= N_{0}'(y)N_{0}(y)\equiv 0,
\]
which identically vanishes since $K$ is normal form, which implies that $N_{0}(y)\equiv 0$. We will prove by induction that $N_{k}'(y)\equiv 0$ for all $k$. Assuming that this is true for all $N_{1}'(y),\cdots,N_{k}'(y)$. Then,
\[
(k+1)\cdot N_{k+1}'(y)=N_{k}'N_{0}+N_{k-1}'N_{1}+\cdots+N_{1}'N_{k-1}+N_{0}'N_{k}.
\]
Thus by the induction hypothesis, $N_{k+1}'(y)\equiv 0$. As a result, we have $N_{k}(y)=c_{k}$ for some constants $c_{k}$. This implies that 
\[
N(x,y)=\sum_{k=1}^{\infty}c_{k}x^{k}.
\]
However, $K$ is a normal form, thus $N(x,0)\equiv 0$. Effectively, $c_{k}=0$ for all $k$, and so $N(x,y)\equiv 0$.

Having proven $N\equiv 0$, we obtain from the first equation $M_{x}\equiv 0$ so that $M(x,y)=M_{1}(y)$ for some function $M_{1}$. Since $K$ is normal, $M(0,y)\equiv 0$, and so $M_{1}(y)\equiv 0$. We have shown that $y''\equiv 0$, and the proof is completed. 
\end{proof}



\bigskip

\noindent{\scriptsize
{\sc Wei Guo {\sc Foo}. Institute of Mathematics, Academia Sinica, 
6F, Astronomy-Mathematics Building, No. 1, Sec. 4, 
Roosevelt Road, Taipei 106319, Taiwan.\\
Email address:} 
\texttt{fooweiguo@hotmail.com}}\\[-10pt]

\noindent {\scriptsize
{\sc Julien {\sc Heyd}. Laboratoire de Math\'{e}matiques d'Orsay, 
UMR 8628 du CNRS, Universit\'{e} Paris-Saclay,
91405 Orsay Cedex, France.\\
Email address:} 
\texttt{julien.heyd@universite-paris-saclay.fr}}\\[-10pt]

\noindent {\scriptsize
{\sc Jo\"el {\sc Merker}. Laboratoire de Math\'{e}matiques d'Orsay, 
UMR 8628 du CNRS, Universit\'{e} Paris-Saclay,
91405 Orsay Cedex, France.\\
Email address:} 
\texttt{joel.merker@universite-paris-saclay.fr}}\\[-10pt]

\vfill
\begin{thebibliography}{225}

{\bf\bibitem{Cartan-1}
{\sc Cartan}}, \'{E}.: 
{\em Sur la g\'{e}om\'{e}trie pseudo-conforme des hypersurfaces de
l’espace de deux variables complexes, I}. 
Ann. Math. Pura Appl. {\bf 11} (1932), 17--90.

\smallskip

{\bf\bibitem{Chern-Moser-1974}
{\sc Chern}}, S. S.; {\sc Moser}, J. K.: 
{\em Real hypersurfaces in complex manifolds}. 
Acta Math. {\bf 133} (1974), 219--271. 
 
\smallskip

{\bf\bibitem{Godlinski-2008}
{\sc Godlinski}}, M:
{\em Geometry of third-order ordinary differential equations
and its applications in General Relativity},
Ph.D. Thesis, arXiv: 0810.2234.

\smallskip

{\bf\bibitem{Godlinski-Nurowski-2009}
{\sc Godlinski}}, M; {\sc Nurowski}, P.:
{\em Geometry of third-order ODEs}, arXiv: 0902.4129.

\smallskip

{\bf\bibitem{Hsu-Kamran-1989}
{\sc Hsu}}, L.; Kamran, N. 
{\em Classification of second-order ordinary differential equations 
admitting Lie groups of fibre-preserving point symmetries},
Proc. London Math. Soc. (3) {\bf 58} (1989), no. 2, 387–416. 

\smallskip

{\bf\bibitem{Merker-2008}
{\sc Merker}}, J.:
{\em Lie symmetries and CR geometry}, J.
Mathematical Sciences, {\bf 154} (2008), no.~6, 817--922.

\smallskip

{\bf\bibitem{Merker-Lie-2015} 
{\sc Merker}},~J. (Editor); {\sc Lie},~S. (Author):
{\em Theory of Transformation Groups I. 
General Properties of Continuous Transformation Groups. 
A Contemporary Approach and Translation}, 
Springer-Verlag, Berlin, Heidelberg, 2015, xv+643~pp.

\smallskip

{\bf\bibitem{Nurowski-Sparling-2003}
{\sc Nurowski}}, P.; {\sc Sparling}, G.A.:
{\em Three-dimensional Cauchy–Riemann structures and second-order 
ordinary differential equations}, 
Classical and Quantum Gravity, Volume 20, Number 23.

\smallskip

{\bf\bibitem{Olver-1995} 
{\sc Olver}}, P.J.:
{\em Equivalence, Invariance and Symmetry}. 
Cambridge, Cambridge University Press, 1995, xvi+525~pp.

\smallskip

{\bf\bibitem{Olver-Pohjanpelto-2008}
{\sc Olver}}, P.J.; {\sc Pohjanpelto}, J.: 
{\em Moving frames for Lie pseudo-groups. }  
Canad. J. Math. {\bf 60} (2008), no.~6, 1336--1386.

\smallskip

{\bf\bibitem{Olver-2018}
{\sc Olver}}, P.J.: 
{\em Normal forms for submanifolds under group actions. 
Symmetries, differential equations and applications}, 
1--25, Springer Proc. Math. Stat., 266, Springer, Cham, 2018. 

\smallskip

{\bf\bibitem{Olver-Valiquette-2018}
{\sc Olver}}, P.J.; {\sc Valiquette}, F.: 
{\em Recursive moving frames for Lie pseudo-groups}, 
Results Math. {\bf 73} (2018), no.~2, Paper No.~57, 64~pp. 

\smallskip

{\bf\bibitem{Ottazzi-Schmalz-2017} 
{\sc Ottazzi}},~A., {\sc Schmalz},~G.:
{\em Normal Forms of Para-CR Hypersurfaces},
Differential Geometry and its Applications
{\bf 52} (2017), 78--93.

\smallskip

{\bf\bibitem{Tresse-1894} 
{\sc Tresse}}, A.:
{\em Sur les invariants diff\'erentiels des groupes
continus de transformations}. Acta Math. 18: 1-88 (1894). DOI: 10.1007/BF02418270.
 
 
{\bf\bibitem{Tresse-1896} 
{\sc Tresse}}, A.:
{\em D\'etermination des invariants ponctuels de l'\'equation 
diff\'erentielle ordinaire du second ordre $y'' = \omega(x,y,y')$},
Hirzel, Leipzig, 1896. 
 
\end{thebibliography}
\end{document}